%% file: 2004-29.tex
\documentclass{gtart_h}  

\input gtoutput

\lognumber{329}

\volumenumber{8}\papernumber{29}\volumeyear{2004}
\pagenumbers{1043}{1078}
\received{5 June 2003}
\published{18 August 2004}
\accepted{11 July 2004}
\proposed{Steven Ferry}
\seconded{Benson Farb, Ralph Cohen}

\usepackage{amsmath,amscd,amssymb}

\newtheorem{thm}{Theorem}[section]
\newtheorem{cor}[thm]{Corollary}
\newtheorem{lemma}[thm]{Lemma}
\newtheorem{corollary}[thm]{Corollary}
\newtheorem{prop}[thm]{Proposition}
\newtheorem*{claim}{Claim}
\newtheorem*{claima}{Claim A}
\newtheorem*{claimb}{Claim B}

\theoremstyle{definition}
\newtheorem*{slc}{The SubLagrangian Construction}
\newtheorem*{cce}{Construction of Characteristic Elements ${\bf v_o, v_1}$}
\newtheorem{definition}[thm]{Definition}
\newtheorem{remark}[thm]{Remark}
\newtheorem{example}[thm]{Example}

\newcommand{\ltt}{{\widetilde{\mcL}}(\zt,2)}

\newcommand{\qlf}{quadratic linking form }
\newcommand{\qnz}{Q^{(n)}_0}
\newcommand{\qnone}{Q^{(n)}_1}
\newcommand{\vnz}{v^{(n)}_0}
\newcommand{\vnone}{v^{(n)}_1}

\newcommand{\Q}{{\mathbb Q}}
\newcommand{\R}{{\mathbb R}}

\newcommand{\N}{{\mathbb N}}
\newcommand{\Z}{{\mathbb Z}}

\newcommand{\eps}{{\varepsilon}}

\newcommand{\mcL}{{\mathcal{L}}}

\newcommand{\mcN}{{\mathcal{N}}}

\newcommand{\mcP}{{\mathcal{P}}}

\newcommand{\mcV}{{\mathcal{V}}}
\newcommand{\bfF}{{\mathbf{F}}}
\newcommand{\bfM}{{\mathbf{M}}}

\newcommand{\ztt}{{\mathbb Z}_2[t]}
\newcommand{\zt}{{\mathbb Z}[t]}

\DeclareMathOperator{\ad}{{Ad}}
\DeclareMathOperator{\ann}{{Ann}}
\DeclareMathOperator{\Arf}{{Arf}}
\DeclareMathOperator{\arf}{{arf}}
\DeclareMathOperator{\cok}{{cok}}

\DeclareMathOperator{\ext}{{Ext}}
\DeclareMathOperator{\ev}{{ev}}
\DeclareMathOperator{\Hom}{{Hom}}
\DeclareMathOperator{\id}{{Id}}
\DeclareMathOperator{\im}{{im}}

\DeclareMathOperator{\length}{{Length}}
\DeclareMathOperator{\od}{{od}}
\DeclareMathOperator{\rank}{{rank}}
\DeclareMathOperator{\unil}{{UNil}}

\begin{document}
\title{The surgery obstruction groups of the\\infinite dihedral
group}
\author{Francis X Connolly\\James F Davis}
\address{Department of Mathematics, University of Notre Dame 
\\ Notre Dame, IN 46556, USA}
\secondaddress{Department of Mathematics, Indiana University
\\ Bloomington, IN 47405, USA}
\asciiaddress{Department of Mathematics, 
University of Notre Dame\\Notre Dame, IN 46556,
USA\\and\\Department of Mathematics, Indiana 
University\\Bloomington, IN 47405, USA}
 
\gtemail{\mailto{connolly.1@nd.edu}{\rm\qua 
and\qua}\mailto{jfdavis@indiana.edu}} 
\asciiemail{connolly.1@nd.edu, jfdavis@indiana.edu} 

\begin{abstract}
This paper computes the quadratic Witt groups (the Wall $L$--groups) of
the polynomial ring $\Z[t]$ and the integral group ring of the infinite
dihedral group, with various involutions. We show that some of these
groups are infinite direct sums of cyclic groups of order 2 and 4. The
techniques used are quadratic linking forms over $\Z[t]$ and Arf
invariants.
\end{abstract}

\asciiabstract{%
This paper computes the quadratic Witt groups (the Wall L-groups) of
the polynomial ring Z[t] and the integral group ring of the infinite
dihedral group, with various involutions. We show that some of these
groups are infinite direct sums of cyclic groups of order 2 and 4. The
techniques used are quadratic linking forms over Z[t] and Arf
invariants.}


\maketitle

\section{Introduction and statement of results}\label{intro}

In this paper we complete the computation of the Wall surgery
obstruction groups for the infinite dihedral group, the $L$--theory of
the polynomial ring $\Z[t]$,  the $L$--theory of the Laurent
polynomial ring
$L_n(\Z[t,t^{-1}])$, with either the trivial involution or the
involution
$t
\mapsto -t$, and the Cappell unitary nilpotent groups for the ring
$\Z$.   The problem of computing these groups is thirty
years old.  We take an historical approach in this introduction which
sets the stage and indicates the interrelation
between the various groups, but has  the drawback of postponing the
discussion of the main results of this paper.  The main results are
\ref{calc}, \ref{expfour}, \ref{computation}  and  \ref{versch thm}.

     Our algebraic computations  are motivated by the following geometric
question: is a homotopy equivalence
     $$
     f\co M \to X_1 \# X_2
     $$
     from a closed $n$--manifold  to a connected sum of two
others {\em splittable}? That is to say, is $M$
expressible as a connected sum $M=M_1\#M_2$
so that $f$ homotopic to a map  of the form
\[
f_1\#f_2\co M_1\#M_2\to X_1\#X_2
\]
where each $f_i\co M_i\to X_i$ is a homotopy equivalence? In
particular, is $M$ itself a connected sum?

Let's restrict now to the case where both $X_i$ are connected
and have cyclic fundamental group of order two.  Cappell
\cite{Cs,CaM} defined an element $s(f)$ in  a 4--periodic abelian
group $\unil_{n+1}(\Z;\Z^\pm,\Z^\pm)$.  The $\pm$'s depend on
the orientability of   $X_1$ and $X_2$ and are often omitted if
both are orientable.  If
$s(f) \not = 0$, then $f$ is not splittable.  If  $s(f) =0 $ and
$n \geq 4$, then $f$ is splittable topologically; if $s(f) =
0$ and $n \geq 5$, then $f$  splits smoothly.  (In the smooth case
one allows connected sum along a homotopy sphere.)  Fixing
$X_1$ and $X_2$ and given $s \in
\unil_{n+1}(\Z;\Z^\pm,\Z^\pm)$, there is a realization result:  there
is a homotopy equivalence $f \co M \to X_1 \# X_2$
with
$s(f) = s$, with $M$ a topological manifold when $n \geq 4$, a smooth
manifold when $n \geq 5$.  A particularly
interesting example is to take $X_1 = X_2 = \R P^4$ and  $0 \neq s
\in\unil_5(\Z;\Z^-,\Z^-)$;
then realization gives a nonsplittable homotopy equivalence $M \to \R
P^4 \# \R P^4$.  The geometric properties of $M$
seem both unexplored and of some interest.

The unitary nilpotent groups $\unil_n(R;A_1,A_2)$ are  defined
for a ring $R$ with involution and $R$--bimodules $A_1,A_2$ with
involution.  We won't need the definitions of these groups, only
their relation to $L$--groups discussed below and the
isomorphisms
$$\unil_n(R;A_2,A_1) \cong \unil_n(R;A_1,A_2) \cong
\unil_{n+2}(R;A_1^-,A_2^-).
$$  Here $A_i^-$ is the bimodule $A_i$, but  with the involution
$a \mapsto - \overline a$.

Associated to a ring with involution  are the  algebraic
$L$--groups $L_n(R)$.  They are 4--periodic.  The definition of
$L_0(R)$ (and
$L_2(R)$) are reviewed in Section \ref{defs};  they are Witt
groups of (skew)-quadratic forms.  $L_{2k+1}(R)$ is the
abelianization of the stable automorphism group of the
$(-1)^k$--hyperbolic form (any form admitting a
Lagrangian).\footnote{These are the so-called $L^h$--groups
measuring the
obstruction to doing surgery up to homotopy equivalence.   They
are defined as in Wall \cite[Chapters 6 and 7]{Wa2}, except one
deletes the requirement that the torsions are trivial.  A
definition of these groups is given in \cite{RI} where they were
denoted
$V_n(R)$. Ranicki \cite{RAT} later gave a definition of
$L_n(R)$ as cobordism classes of $n$--dimensional
quadratic Poincar\'e complexes over $R$.}

     Fix a group $G$ and a homomorphism $w \co G \to \{\pm 1 \}$. There is an
induced involution
$\sum a_g g \mapsto \sum a_g w(g) g^{-1}$on
$\Z G$.  The associated groups
$L_n(\Z G,w)$ are key ingredients in the classification
of closed, oriented manifolds with fundamental group $G$ and
orientation character $w$.

Parallel to the work of Stallings \cite{St} and Waldhausen
\cite{Wald} in algebraic $K$--theory,
Cappell
\cite{Cu,Cal}  studied the
$L$--groups of amalgamated free products and showed that if $H$ is a
subgroup of groups $G_1$ and $G_2$, then $\unil_n(\Z H; \Z [G_1 - H],
\Z [G_2 -H])$ is a summand of $L_n(\Z [G_1
*_H G_2])$, and that the $L$--group modulo the UNil--term fits into a
Mayer--Vietoris exact sequence.  Farrell \cite{Fa}
showed that the UNil--term has exponent at most four. However he was
unable to find an element $\alpha \in \unil_n(R;A,B)$
for which $2\alpha \not = 0$.   Cappell proved that the
$\unil$--term vanishes provided that the inclusions
$H
\hookrightarrow G_i$ are {\em square root closed}, ie,
$g \in G_i$ and $g^2 \in H$ implies $g \in H$.

The infinite dihedral group is
$$
D_\infty =
\Z_2 * \Z_2 = \langle a_1,a_2~|~a_1^2 =e=a_2^2 \rangle
= \langle g,t ~|~ t^2 = e, tgt^{-1} = g^{-1}\rangle.
$$
     Let $w \co D_\infty \to \{\pm 1\}$ be a homomorphism.  The
$L$--groups
$L_n(\Z[D_\infty],w)$ and the corresponding
$\unil$--groups
$\unil_n(\Z;\Z^{w(a_1)},\Z^{w(a_2)})$  seem particularly
fundamental.  First the infinite dihedral group is  the
simplest group which is not square root closed.  Second, due to
the isomorphism conjecture of Farrell--Jones \cite{FJ}
(generalizing the Borel--Novikov conjectures of manifold theory),
attention has been recently focused on the infinite dihedral
group.  The isomorphism conjecture roughly states that $L_n(\Z
\Gamma)$ depends on the $K$-- and $L$--theory of virtually cyclic
subgroups and homological data depending on $\Gamma$.  A group
$G$ is virtually cyclic if either $G$ is finite, or $G$ surjects
onto $\Z$ with finite kernel, or $G$ surjects onto $D_\infty$
with finite kernel.  The $L$--theory in the first two cases has
been examined in detail  \cite{HT,S,RII,RIII}.  Therefore
$L_n(D_\infty)$ is the next obvious case to consider.

In this paper we are writing a conclusion to the long tale of the
computation of $\unil_n(\Z;\Z^{\pm},\Z^{\pm})$.
      Cappell showed
           $\unil_2(\Z;\Z,\Z)$ was infinitely generated \cite{C} and
announced that
$\unil_0(\Z;\Z,\Z)=0$ \cite{Cu}.   Connolly--Ko\'zniewski
\cite{CK} obtained an isomorphism $\unil_2(\Z;\Z,\Z) \cong (\Z_2)^\infty$
and also showed that\break  $ \unil_0(\Z;\Z,\Z)=0.$   Connolly--Ranicki \cite{CR}
showed $\unil_1(\Z;\Z,\Z) = 0$ and computed $\unil_3(\Z;\Z,\Z)$ up to
extension,
and thereby showed that it was infinitely generated.  Andrew Ranicki  
outlined the
construction and detection of  an element of $\unil_3(\Z;\Z,\Z)$  of
exponent 4 in a letter  \cite{Rl} to the first author.  After a
preliminary version of this paper was circulated, Banagl--Ranicki
\cite{BR} gave an independent complete computation of
$\unil_3(\Z;\Z,\Z)$ using generalized Arf invariants.

In this paper we give complete  computations for all
$n$ as well as doing the non-oriented
case.

But before we discuss our computations we pause and explain how
computations of the unitary nilpotent group give explicit
computations of the $L$--theory of the infinite dihedral group.
We rely on the Mayer--Vietoris exact sequence (see Cappell \cite{Cu}):
\begin{multline*}
      \dots \to L_n(\Z)\to L_n(\Z[\Z_2],w_1)\oplus  L_n(\Z[\Z_2], w_2)\\ \to
\frac{L_n(\Z[D_\infty], w)}{ \unil_n(\Z;\Z^{\eps_1},\Z^{\eps_2})}\to
L_{n-1}(\Z)\to \cdots
\end{multline*}
where $w_i = w|_{\langle a_i \rangle}$
     and $ \eps_i=w(a_i)=\pm1$. We assume $\eps_2 \geq \eps_1$ and write
$\widetilde{L}_n(\Z[G],w)$ for the cokernel of the natural map $ L_n(\Z)\to
L_n(\Z[G],w).$ The above sequence, and the calculational results in
Wall \cite[Chapter 13A]{Wa2},
     quickly lead us to the following equations.
\[
L_n(\Z[D_\infty], w)=\widetilde{L}_n(\Z[\Z_2],w_1)\oplus L_n(\Z[\Z_2],w_2)\oplus
\unil{_n(\Z; \Z^{\eps_1},\Z^{\eps_2})}
\]
unless $n \equiv 1 \text{ mod } 4$, and $w_1, w_2$ are both nontrivial; in this
case, we get:
\[
L_1(\Z [D_\infty], w)=
\unil{_1(\Z; \Z^-, \Z^-)}\oplus L_0(\Z).
\]
For the values of $L_n(\Z[\Z_2],w)$  see Wall \cite[Chapter 13A]{Wa2}.

        There is another relation between the unitary nilpotent groups
and $L$--groups which will be crucial to our
computations.  Let
$R[t^{\pm}]$ denote the polynomial ring
$R[t]$ with
$+$ involution
$\sum r_i t^i
\mapsto
\sum
\overline{r_i} t^i$  or the $-$ involution $\sum r_i t^i \mapsto \sum
(-1)^i \overline{r_i} t^i$.  Let
$\varepsilon_0 \co R[t] \to R$ be the map of rings given by
$\varepsilon_0(f(t)) = f(0)$; it is a split surjection with the
splitting   $s_0\co R \to R[t]$ given by constant
polynomials.  Define $$N^{\pm }L_n(R) =  \ker (\varepsilon_0 \co
L_n(R[t^{\pm }]) \to L_n(R)).$$ Then $L_n(R[t^{\pm }]) = L_n(R) \oplus
N^{\pm }L_n(R)$.  Note that $L_n(\Z) \cong \Z,0,\Z_2,0$ for $n \equiv
0,1,2,3
\pmod 4$ so our computation of $N^{\pm }L_n(\Z)$  given below also
computes  $L_n(\Z[t^{\pm}])$.
\begin{thm}[Connolly--Ranicki]\label{unil=nl}
There is an isomorphism, natural in $R$,
$$
r\co \unil_n(R;R^{\pm },R) \to N^{\pm }L_n(R).
$$
\end{thm}
Connolly--Ranicki \cite[Theorem A]{CR} prove this in the $+$ case; we remark here
that their formula \cite[Definition 2.13]{CR}  for $r$ and their proof that $r$
is an isomorphism apply equally well in the $-$ case.  Note that
this implies a 2--fold periodicity
$$
N^-L_n(R) \cong \unil_n(R;R^-,R) \cong \unil_{n+2}(R;R,R^-)  \cong
N^-L_{n+2}(R).
$$
These $NL$--groups are analogous to Bass' nilpotent $K$--groups
$NK_i(R) =\linebreak \ker(\eps\co K_i(R[t] \to K_i(R))$ occurring in the
fundamental theorem of algebraic $K$--theory \cite{B}.

The following theorem, which is an easy consequence of Theorem \ref{3n3}, provides the calculation of $N^\pm L_{2k}(\Z)$.

\begin{thm} \label{even_comp}  There are isomorphisms of abelian groups:
\begin{enumerate}
\item  $\varepsilon_{0*} \co L_0(\zt) \xrightarrow{\cong} L_0(\Z)$.  Thus $NL_0(\Z) = 0$.
\item  $ P_2^{\zt} \co t\ztt/\{p^2 - p : p \in t\ztt\}  \xrightarrow{\cong} NL_2(\Z) $.
\item  $P_2^{\Z[t^-]} \co t^2\Z_2[t^2]/\{p^2 - p : p \in t^2\Z_2[t^2]\}
\xrightarrow{\cong}  N^-L_2(\Z).$
\item  $ Q_0^{\Z[t^-]}\co  t^2\Z_2[t^2]/\{p^2 - p : p \in t^2\Z_2[t^2]\}
\xrightarrow{\cong} N^-L_0(\Z)$.

\end{enumerate}
\end{thm}

The maps $P_2^{\Z[t^\pm]}$ and $Q_0^{\Z[t^-]}$ are defined in Section \ref{computesection}.  The inverse maps in (2), (3), and (4) are all essentially given by the Arf invariant  of
the function field $\Z_2(t)$.  Also in Section \ref{computesection}, we compute $L_2(\Z_2[t]) = L_0(\Z_2[t])$ and show it is isomorphic to $L_2(\Z[t])$.

But why do these polynomials $p(t)$ appear  in the computation of the
$L$--groups?   None of the
     groups above are $\Z[t]$--modules, but the above  isomorphisms seem to
be more than isomorphisms of abelian groups.  We explain that now.

For each integer $i>0$, and each ring $R$, we have a ring endomorphism:
\[
V_i\co R[t]\to R[t]\quad \qquad V_i(p(t)) = p(t^i).
\]
Note  $\eps_0\circ V_i = \eps_0$.
The resulting monoid of endomorphisms of $L_n(\Z[t])$,
\[
\mathcal{M}= \{V_1, V_2, V_3,\dots \} ; \quad V_iV_j = V_{ij}
\]
       therefore  makes $L_n(\Z[t])$ and $NL_n(\Z)$ 
modules over  the
\emph{Verschiebung  Algebra},
\[
\mcV =\Z[\mathcal{M}],
\]
a polynomial ring on $\{V_p: p \text{ is a
prime}\}$:  $\mcV=\Z[V_2, V_3, V_5,\dots] $ .
       The subalgebra indexed by the odd primes, $\mcV_{\od}:=\Z[V_3,
V_5, V_7, V_{11},\dots]$  acts similarly on
$N^-L_n(\Z)$.   (Note that for $i>0$ even,  $V_i \co  \Z[t^-] \to \Z[t^-]$ is not a map of rings with involution.)  The map in Theorem \ref{even_comp}(2) is a map of $\mcV$--modules and the maps in Theorem  \ref{even_comp}(3) and \ref{even_comp}(4) are maps of $\mcV_{\od}$--modules.
We have the following reformulation of Theorem \ref{even_comp}:
\begin{thm}  \begin{align*}
&NL_{2k}(\Z_2) \cong NL_{2}(\Z) \cong \mcV/\langle 2, V_2 - 1\rangle \\
&N^-L_{2}(\Z) \cong N^-L_{0}(\Z) \cong \mcV_{\od}/\langle 2\rangle \\
&NL_0(\Z) = 0
\end{align*}
\end{thm}
To attack the odd-dimensional $L$-- and $\unil$--groups, we use the
classical technique of quadratic linking forms.  Given a
ring with involution $R$ and a central multiplicative subset $S =
\overline S$ of non-zero divisors, one can define
$L_{2k}(R,S)$ to be the Witt group of $(-1)^k$--quadratic linking
forms on finitely generated $S$--torsion $R$--modules of
length one (see \cite{RES}).  Furthermore, one can identify
$L_{2k}(R,S)$ with the relative $L$--group $L_{2k}(R \to
S^{-1}R)$.  There is an analogous theory for $L_{2k+1}$.

For a ring with involution $R$ where $2$ is not a divisor of zero, define
$\langle 2 \rangle = \{2^i : i \geq 0\} \subset R$ and
$$
N^{\pm}L_{n}(R,\langle 2 \rangle) = \ker \eps_0 \co  L_n(R[t^{\pm}]),
\langle 2 \rangle) \to L_n(R, \langle 2 \rangle).
$$

\begin{prop}  \label{local} For any Dedekind domain with involution $R$ where $2$ is not a
divisor of zero,
$$
N^{\pm}L_{n}(R,\langle 2 \rangle) \cong N^\pm L_{n-1}(R).
$$
\end{prop}

\begin{proof}   By comparing the long exact sequence of the ring map
$$R[t^\pm] \to R[1/2][t^\pm]$$ with that of $R \to R[1/2]$, one
obtains a localization  exact sequence 
$$
\cdots \to N^{\pm}L_n(R[1/2]) \to N^{\pm}L_n(R,\langle 2 \rangle) \to
N^{\pm}L_{n-1}(R) \to N^{\pm}L_{n-1}(R[1/2]) \cdots
$$
Since 2 is a unit, $N^\pm L_k(R[1/2]) = N^\pm L^k(R[1/2])$.  
Connolly--Ranicki  show that for a Dedekind domain $\Lambda$ with involution, $N^\pm L^k(\Lambda) = 0$
(see \cite[Proposition 2.11, Proposition 2.19, and the discussion after Proposition 2.19]{CR}). 
\end{proof}

These Witt groups of quadratic linking forms are the main object of
study in our paper.  They occur so often that introduce new notation
for them.
\begin{align*}
\mcL(\Z[t^\pm],\langle 2 \rangle ) &= L_0(\Z[t^\pm],\langle 2 \rangle )\\
\widetilde \mcL(\Z[t^\pm],\langle 2 \rangle ) &= N^\pm L_0(\Z,\langle
2 \rangle )\\
\mcL(\Z,\langle 2 \rangle ) &= L_0(\Z,\langle 2 \rangle ) \cong \Z_2
\oplus \Z_8
\text{ (see \cite[Theorem 5.2.2]{Sc})}
\end{align*}

We also use the notation $\mcL(\Z[t^\pm], 2^n )$ to denote the Witt
group of quadratic linking forms on length one $\Z[t]$--modules of
exponent $2^n$ (see Section \ref{defs}.)  To study these groups 
we use the method of characteristic elements.  

A key technical result for
us in the following devissage result (see Section \ref{S3}).

\begin{thm}  \label{devissage}
$\widetilde \mcL(\Z[t^\pm],2) \cong
\widetilde
\mcL(\Z[t^\pm],\langle 2\rangle).$
\end{thm}

We next present some examples of quadratic linking forms on
$\Z[t]$--modules with exponent 2.

\begin{definition}\label{intronpg} For polynomials $p, g \in \zt$, define
the quadratic linking form
\begin{equation}\label{npgdef}
\mcN_{p,g}  =
\left(
\ztt^2,
\begin{pmatrix}
p/2 & 1/2\\
1/2& 0
\end{pmatrix},
\begin{pmatrix}
          p/2\\
          g
\end{pmatrix}
\right)
.
\end{equation}

By this notation we mean that if $\{e_1, e_2\}$ is the standard basis
of $\ztt ^2$ over $\ztt$, then the $2\times 2$ matrix
is
$(b(e_i,e_j))$ and the column vector is $(q(e_i))$.  If $p(0) \in
4\Z$, or $g(0)\in 2\Z$, then $[\mcN_{p,g}]$ is an element of
$\widetilde
\mcL(\zt,\langle 2\rangle)$.
\end{definition}

We wish compute these Witt groups as modules over the Verschiebung
Algebra.  To this end, note  $t\Z_2[t]$ is a $\mcV$--module in
the  obvious way: $V_n(p(t)) = p(t^n)$. But, as a less obvious module,
$[t\Z_2[t]]$ shall denote the abelian group $t\Z_2[t]$, equipped with the
following $\mcV$--module structure:
\begin{align}\label{vstructure}
V_{2n+1}(p(t)) &= p(t^{2n+1});\\
V_{2n}(p(t)) &= 0,\notag
\end{align}
for all $n\geq 0$, and all $p\in t\ztt$.

We are now is a position to state our  main theorem. (See
Section \ref{section-answer} for the proof.)

\begin{thm}  \label{computation}
\begin{enumerate}

\item
There is an isomorphism of $\mcV$--modules,
\[
j_1 + j_2 \co  \frac{t\Z_4[t]}{(2V_2-2)t\Z_4[t]}
\oplus [ t\Z_2[t]] \to \widetilde \mcL (\zt,\langle 2\rangle)
\]
given by $j_1[tp] = [\mcN_{tp,1}]$ and $j_2[tp] = [\mcN_{1, tp}]- [\mcN_{t, p}]$.
\item There is an isomorphism of $\mcV_{\od}$--modules,
\[
j\co t^2\Z_2[t^2]\to \widetilde{\mcL}(\Z[t^-], \langle 2\rangle); \qquad j(t^2p(t^2))
=[\mcN_{t^2p(t^2), \,1}].
\]
\end{enumerate}
\end{thm}

\begin{cor}\label{versch thm}
\[
 \widetilde \mcL (\zt,\langle 2\rangle) =\bigoplus_{k=-1}^\infty \mcV\cdot b_k;
\qquad \widetilde{\mcL}(\Z[t^-], \langle 2\rangle)= \bigoplus_{j=1}^\infty \mcV_{\od}\cdot c_j
\]
where $b_{-1}=[ \mcN_{t,\;1}]   $ and $b_k = [\mcN_{1,\; t^{2^k}}] - 
[\mcN_{t, \;t^{2^k-1}}]$ for all $k\geq 0$, and $c_j= \mcN_{t^{2^j}, 1}$, 
for all
$j\geq 1$.  The annihilator ideals of these elements are:
\[
\ann(b_{-1}) =\langle 4,2V_2-2\rangle ;\; \ann(b_k)= \langle 2,V_2
\rangle  \text{ for } k\geq
0;\; \ann_{\,\mcV_{\od}}(c_j)= \langle 2 \rangle  .
\]
\end{cor}
\begin{proof} Note that $t\Z[t]$ is a free rank one
$\mcV$--module.  In $ \frac{t\Z_4[t]}{(2V_2-2)t\Z_4[t]} $ ,  a $\mcV$--generator
is $t$, with annihilator ideal $\langle 4, 2V_2-2\rangle $. Also
$[t\ztt]$ is a free
module over $\mcV/\langle 2,V_2\rangle $,  with basis $\{t^{2^k}: k=
0,1,2, \dots\}$.
Finally the set $\{t^{2^j}: j=1,2,3, \dots\}$ is a $\mcV_{\od}/2\mcV_{\od}$
basis for
$t^2\Z[t^2]$. The result now follows from Theorem \ref{computation}.
\end{proof}

As a consequence of our computation we have the following
corollary (cf \cite{Rl}). 
\begin{corollary}\label{expfour}
There is an element $\alpha \in \unil_3(\Z;\Z,\Z)$ of order 4.
\end{corollary}

We give a separate, elementary proof of this using Gauss sums in Section \ref{expfoursection}.

The heart of our paper is  the following set of calculations
which we summarize as follows:
\begin{thm}\label{calc}
\begin{align*}
     \unil_0(\Z;\Z^-,\Z)&=
\unil{_2(\Z;\Z^-,\Z)}\cong\mcV_{\od}/\langle 2 \rangle
\quad
\text{(Theorem \ref{3n3})}
\\
     \unil_1(\Z;\Z^-,\Z)&=
 \unil{_3(\Z;\Z^-,\Z)}\cong\bigoplus_{i=1}^\infty\mcV_{\od}/\langle 
2\rangle
\qquad\text{(Corollary
\ref{versch thm})}
\\  \unil{_0(\Z;\,\Z \,,\Z)}&=  0 \qquad\text{(Theorem
\ref{3n3})}
\\
\unil{_1(\Z;\,\Z \,,\Z)}&= 0 \qquad\text{(Proposition \ref{un1})}
\\
\unil{_2(\Z;\,\Z \,,\Z)}&\cong\; \mcV/\langle 2,V_2-1\rangle \qquad \text{(Theorem
\ref{3n3})}\\
     \unil{_3(\Z;\,\Z \,,\Z)}&\cong \; \mcV/\langle 4,
2V_2-2\rangle \oplus\bigoplus_{i=0}^\infty\mcV/\langle 2,V_2\rangle
\qquad\text{(Corollary
\ref{versch thm})}
\end{align*}
\end{thm}
Finally, let $\Z[t,t^{-1}]^\pm$ be the Laurent polynomial ring
$\Z[t,t^{-1}]$ with the involution $t \mapsto \pm t$.  The following
formulas and the results of this paper   compute
$L_n(\Z[t,t^{-1}]^\pm)$.
\begin{align*}
L_n(\Z[t,t^{-1}]) & \cong L_n(\Z)^2 \oplus NL_n(\Z)^2 \\
L_n(\Z[t,t^{-1}]^-) & \cong L_n(\Z)\oplus L_{n+2}(\Z) \oplus N^-L_n(\Z)^2
\end{align*}
The first formula is proved by Ranicki in \cite{RIV} in the oriented case, and
the second can be proven with similar techniques.

\section{Definitions} \label{defs}

Here we
present a unified framework including both the Witt group of quadratic
linking forms over a ring  $R$, and the surgery obstruction groups
$L_{2k}(R)$.

Let $R$ be a ring. An $R$--module $M$ has {\em length one} if  there is a
short exact sequence $0
\to F_1 \to F_0 \to M \to 0$, where $F_0$ and $F_1$ are finitely
generated free $R$--modules.  A submodule $L$ of a length one module $M$ has
{\em colength} one if
$M/L$ has length one. This implies $L$ has length one.

         A {\em ring with involution} is a ring $R$ with a function $-\co  R
\to R$ satisfying $\overline{1} = 1$, $\overline{\overline{r}}=r$,
$\overline{r + s} = r + s$, and $\overline{rs} = \overline{s}
~\overline{r}$ for all $r,s \in R$.
An {\em $R$--bimodule
with involution} is
an $R$--bimodule $A$ with $\Z$--automorphism $- \co  A \to A$ of order
2 satisfying $\overline{r a s} = \overline{s}~
\overline{a} ~ \overline{r}$ for all $r,s \in R$, $a \in A$.

Let $(R, A)$ be a ring $R$  with involution, together with some
$R$--bimodule with involution.   If $M$ is a left $R$--module, then
$M^\wedge:=\Hom_R(M,A)$ is also a left $R$--module if we set $(r \varphi)(m) =
\varphi(m) \overline{r}$.
A {\em symmetric form over $(R, A)$}
is  an $R$--module $M$ and a function
$\lambda \co  M \times M \to A$ satisfying:
\begin{align*}
\lambda(x,y) & = \overline{\lambda(y,x)} \\
\lambda(x,y+y') & = \lambda(x,y) + \lambda(x,y') \\
\lambda(x, ry) & = r\lambda(x, y)
\end{align*}
for all $r \in R$, $x, y, y' \in M$.  If the $R$--map
\begin{align*}
\ad \lambda \co  M &\to M^\wedge=\Hom_R(M,A)
\qquad x  \mapsto \;(y \mapsto \lambda(x, y))
\end{align*}
is an isomorphism, then $(M,\lambda)$ is {\em nonsingular}.

A {\em (nonsingular) quadratic form over} $(R, A)$ is a
(nonsingular) symmetric form $(M,\lambda)$ over  $(R, A)$
and a function
$\mu\co M \to A/ \{a- \overline{a}:a \in A\} $ so that
\begin{align*}
\lambda(x, x) &= \{1+*\}\mu(x)\\
[\lambda (x, y)] &= \mu(x+y)
- \mu(x) -
\mu(y)
\\
\mu(rx) &= r \mu(x)\overline{r}
\end{align*}
for all $r \in R, \;x,y\in M$.  Here $[~] \co  A \to A/ \{a-
\overline{a}: a \in A \}, a \mapsto [a]$ is the quotient map and
$\{1+*\}\co  A/
\{a- \overline{a}: a \in A \} \to A$ is given by
$[a]
\mapsto a + \overline{a}$.

       A {\em skew-quadratic }(or {\em$(-1)$--quadratic}) form  over
$(R, A)$ is a quadratic form over $(R, A^-)$.
Here
$A^-$ is the bimodule $A$, but  with the involution:
$a \mapsto -\overline{a}$.

\begin{definition}\label{triple}
      Suppose $R$ is a ring with involution, $A$ is an
$R$--bimodule with involution,
      and $\bfM$ is a class of left  $R$--modules. We assume  $(R,
A,
\bfM)$
      satisfies the following properties:
\begin{enumerate}
\item  If $M\in \bfM$, then $M^\wedge \in \bfM$, and
the double duality map $M\xrightarrow{D}M^{\wedge \wedge}$ is an
isomorphism.
\item The direct sum of two modules in $\bfM$ is in $\bfM$.   $\{0\}$ is in
$\bfM$.
\item If $0\to M_1\to M_2\to M_3\to 0$ is  exact and $M_2, M_3$ are
in $\bfM$, then $M_1$ is also in $\bfM$,  and $0\to M_3^\wedge\to M_2^\wedge\to
M_1^\wedge \to 0$ is also exact.
\end{enumerate}
      A \emph{$(R,A, \bfM)$--form} is a nonsingular quadratic form
$m=(M,\lambda,
\mu)$ over
$(R,  A)$, such that $M\in \bfM$.
\end{definition}

Next we give the examples we
      care about.  It is clear that the first example below satisfies
the conditions of \ref{triple}.  
Proposition \ref{tripl} proves that the second example also satisfies the
conditions.
\begin{itemize}
\item $(R, R, \bfF)$. Here $R$ is a ring with involution; $A=R; \;
\bfF$  is the class of all finitely generated stably free $R$--modules. We 
will also use $(R, R^-, \bfF)$.
\item $(R, \Q R/R, \bfM_{2^n})$. Here $R$ is a ring with involution which is
torsion free as an abelian group and
$A=\Q R/R$, where $\Q R $    denotes the localization of $R$ obtained
by inverting the positive integers in $R$;\; $\bfM_{2^n}$ is the
class of  all length one $R$--modules $M$ for which
$2^n M = 0$.

We will also use $(R, \Q R/R, \bfM\langle 2\rangle),$ where $\bfM
\langle 2\rangle:=
\cup_{n\geq 0} ~\bfM_{2^n}$.
\end{itemize}

A {\em subLagrangian} for a $(R,A, \bfM)$--form $m =(M, \lambda, \mu)$
is a
       submodule $L\subseteq M$ such that $M/L \in \bfM$, $\mu(L)=0$, and $L
\subseteq L^\perp$, where
\[
L^\perp  := \{ x \in M : \lambda(x,L) = 0\}.
\]
If $L=L^\perp$, we say  $L$ is a \emph{Lagrangian  for $m$}, and that \emph{$m$
admits a Lagrangian}.

Define an equivalence relation on the  collection of
$(R,A,\bf M\rm)$--forms by:
$$
m \sim m' \quad \text{if} \quad m \oplus m_0 \cong m' \oplus m_0'
$$
where $m_0$ and $m_0'$ are $(R,A,\bf M\rm)$--forms   which
admit Lagrangians.  (Here $\cong$ means ``is isometric to''.)
\begin{definition}  The \emph{quadratic Witt group $QW(R,A,\bfM)$} is the abelian
group of equivalence classes of  $(R,A,\bf M\rm)$--forms. Addition is
orthogonal  direct sum.
\end{definition}
      The negative of
$[M,\lambda,\mu]$  is $[M,-\lambda,-\mu]$, since their sum admits the diagonal
Lagrangian $\{(x,x): x \in M\}$.

\begin{slc}
Given a  $(R, A,\bfM)$--form $m =(M, \lambda, \mu)$
      and a subLagrangian $S$ for $m$,
we define an induced   quadratic
form
\[
m_S = (S^\perp/S, b_S,\mu_S):  \qquad b_S([x],[y]) = b(x,y) ;\quad \mu_S([x]) =
\mu(x).
\]
\end{slc}

Lemma \ref{sublag} shows that  $m_S$ is again a $(R, A, \bfM)$--form and $[m]=[m_S]\in
QW(R,A,\bfM)$.

The Wall  surgery obstruction groups of $R$ (in even
dimensions) are:
\[
L_0(R)  = QW(R,R,\bfF )\qquad
L_2(R) = QW(R,R^-,\bfF ),
\]
where $\bfF$ is the class of finitely generated, stably free $R$--modules.

We are going  to define the \emph{Witt groups  of quadratic linking
forms} $\mcL(R,2^n)$ and $\mcL(R,\langle 2\rangle)$, so that,
\begin{align*}
\mcL(R,2^n) &\cong QW(R,\Q R/R,\bfM_{2^n})\\
\mcL(R,\langle 2\rangle) &\cong QW(R,\Q R/R,\bfM\langle 2
\rangle)
\end{align*}

       But for historical reasons
(see for example
\cite{Wa1}), we will first change coordinates, and do all of our work in the
classical regime of \emph{quadratic linking forms}.

        \begin{definition}
Assume $R$ is a ring with involution which is torsion
free as an abelian group.

A \emph{quadratic linking form}
over
$R$ is a triple $(M,b,q)$ such that $M\in
\bfM_{2^n}$ for some $n$,  $(M,b)$ is a nonsingular symmetric form over $(R, \Q
R/R)$,  and
$q \co  M
\to
\Q R/(1+*) R$ is a function satisfying:
\begin{align*}
          [q(x)] &= b(x,x) \in \Q R/R ,\\
       q(x+y) - q(x) - q(y) &= \{1+*\}b(x,y) \in \Q R/(1+*)R,\\
       q(rx) &= rq(x)\overline{r}.
\end{align*}
Here 
$$
\{1+*\}\co \frac{\Q R/R}{(1-*)\Q R/R}\to \frac{(1+*)\Q R}{(1+*)R}
$$ 
is
the isomorphism: $[x]\mapsto (x+\overline{x}) \mod{(1+*)R}$.

Note  $2^{n+1}q(x)=0$, whenever $2^nb(x,x)=0$.  Also 
$q(M)\subset
\frac{(1+*)\Q R}{(1+*)R}$ .
\end{definition}
In other words, a \qlf $(M,b,q)$ is any triple which can be written
in the form $(M,b,\{1+*\}\mu)$, where
$(M, b,
\mu)$ is a $(R, \Q R/R, \bfM\langle 2\rangle)$--form.

Let $\Lambda(R,  2^n)$ denote the set of isometry classes of
quadratic linking forms
      with exponent
$2^n$.
The rule \begin{equation}\label{corresp}
(M,\lambda, \mu)\leftrightarrow
(M,\lambda, \{1+*\}\mu)
\end{equation}
gives  a one to one correspondence between the set of isometry classes of\linebreak
$(R, \Q R/R,
\bfM_{2^n})$--forms and $\Lambda(R,  2^n)$.

A Lagrangian for a quadratic linking form
$(M,b,\{1+*\}\mu)$ is defined to be a Lagrangian for the $(R, \Q R/R,
\bfM\langle 2\rangle)$--form
$(M,b,\mu)$.

Two elements $m, m'\in \Lambda(R, 2^n)$ are \emph{equivalent} if
$ m \oplus m_0 \cong m' \oplus m_0'
$
where $m_0$ and $m_0'$ are  elements of  $\Lambda(R, 2^n)$ which admit
Lagrangians.
\[\text{Set: } \qquad\mcL(R,2^n) = \Lambda(R, 2^n)/\sim,
\]
       an abelian  group under orthogonal direct sum.   \eqref{corresp}
induces an isomorphism
\[ QW(R,\Q R/R,\bfM_{2^n}) \cong \mcL(R,2^n).
\]
Similarly, the above equivalence relation on each of the $\Lambda(R,2^n)$'s
defines an equivalence relation on their union,
\[
\Lambda(R,\langle 2\rangle)=\cup_{n=1}^{\infty}\Lambda(R,2^n),
\]
and the abelian group of equivalence classes,
\[
\mcL(R, \langle 2\rangle ):=\Lambda(R,\langle 2\rangle)/\sim
\]
is canonically isomorphic to $QW(R, \Q R/R,
\bfM\langle 2\rangle)$ by the rule \eqref{corresp}.

One could define a corresponding group of skew-quadratic linking forms, but we do not do so here, because we show in Proposition \ref{un1} that when $R = \Z$, any such form admits a Lagrangian.

\section{The proof of Corollary \ref{expfour}}
\label{expfoursection}

In this section  we give a
short  proof, independent of the rest of this
paper, of the fact  that $\widetilde \mcL(\zt,2) \cong
\unil_3(\Z;\Z,\Z) \cong L_3(\Z[t])$ has an element of order 4.
We first review the isomorphism \cite[5.2.2]{Sc}.
\[
\text{Rk} \oplus \text{GS} \co  \mcL(\Z,\langle 2 \rangle)
\xrightarrow{\simeq} \Z_2 \oplus \Z_8.
\]
The rank homomorphism is
\[
\text{Rk}[M,b,q] =
\begin{cases}
0 & \text{if $|M|= 2^{2k}$} \\
1 & \text{if $|M|= 2^{2k+1}$}.
\end{cases}
\]
The Gauss sum homomorphism is
$$
\text{GS}[M,b,q] = [k] \in \Z_8
$$
where
$$
\frac{1}{\sqrt{|M|}} \sum_{x\in M} e^{\pi i q(x)} = e^{2\pi i k}.
$$
Now let $\alpha = \left[
\Z_2[t]^2,
\begin{pmatrix}
1/2 & 0\\
0  & 1/2
\end{pmatrix},
\begin{pmatrix}
          1/2\\
          t-1/2
\end{pmatrix}
\right] \in \mcL(\Z[t],\langle 2 \rangle)$.  Note that
$\varepsilon_0(\alpha)$ has the diagonal Lagrangian $L =
\{(0,0),(1,1)\}$, so
$\alpha \in\widetilde \mcL(\Z[t],\langle 2 \rangle)$.   Consider now
the ring map
\[\varepsilon_1 \co  \Z[t] \to \Z, \quad f(t) \mapsto f(1).
\]
$$
\text{Then:}\qquad  (\text{Rk} \oplus \text{GS})\varepsilon_1(\alpha)
= 0 \oplus
2
\in
\Z_2
\oplus \Z_8,
$$
so $\alpha \in\widetilde \mcL(\Z[t],\langle 2 \rangle)$ is an  element of
order at least 4.  One can show $4\alpha =0$ by  quoting Farrell's
     Exponent Four Theorem \cite{Fa}  or by showing directly that
$4\alpha$ has a Lagrangian.

\section{Computation of $\unil_{2n}(\Z; \Z, \Z),
\unil_{2n}(\Z; \Z^-, \Z)$, and $\mcL^{\ev}(\Z[t^{\pm}], 2)$}
\label{computesection}

According  to Theorem \ref{unil=nl}, $$L_{n}(\Z[t^{\pm}]) = L_n(\Z) \oplus N^{\pm}L_n(\Z) =
L_{n}(\Z)\oplus
\unil_{n}(\Z; \Z^{\pm}, \Z).$$
In this section we compute  $L_{2n}(\Z[t^{\pm}])$, and therefore
$\unil_{2n}(\Z;
\Z^{\pm}, \Z)$. We also compute the group $\mcL^{\ev}(\Z[t^{\pm}],2)$ of
       quadratic linking forms of   exponent $2$ with even type, as defined
below.

Throughout this section, $R$ denotes  a  ring with
involution.

\begin{definition}\label{clod}
Fix an integer $ n>0$. Assume $R$ is torsion free as an abelian group.
Let $\Lambda^{\ev}(R,2^n)$ be the collection of all  quadratic linking
forms
$m=(M, b, q )$ for which $2^nM=0$ and $2^nq(x)=0$  for all $x \in M$.
Quadratic  linking forms in $\Lambda^{\ev}(R, 2^n)$ are called {\em forms of
exponent
$2^n$ with even type}.

We
define an equivalence relation on $\Lambda^{\ev}(R,2^n)$ by:
\[
m\sim m' \Leftrightarrow m\oplus m_0\cong m'\oplus m'_0
\]
for some $m_0, m'_0$ in $\Lambda^{\ev}(R,2^n)$ admitting Lagrangians.
The abelian group of equivalence classes of forms of
exponent
$2^n$ with even type is:
\[
\mcL^{\ev}(R, 2^n) = \Lambda^{\ev}(R, 2^n)/\sim
\]
\end{definition}

Forgetful maps define homomorphisms,
\[
\mcL(R,2^{n-1})\xrightarrow{j_{n-1}}
\mcL^{\ev}(R,2^{n })\xrightarrow{\iota_n }
\mcL(R,2^{n })\xrightarrow{i_n}\mcL(R, \langle 2\rangle).
\]
The map  $\iota_n \circ  j_{n-1}$ is  mentioned  in Theorem
\ref{ses2q} below.

There is a natural map
$$
\alpha_*\co L_0(R/ 2R )\to \mcL^{\ev}(R, 2), \qquad[M,\lambda,\mu]
\mapsto [M,\alpha\circ \lambda,\beta \circ\mu]
$$
where
$$
R/2R \xrightarrow{\alpha} \Q R/R, \quad [x]  \mapsto [x/2]
$$
$$
\frac{R/2R}{ (1-*)(R/2R)}\xrightarrow{\beta}
\Q R/(1+*)R, \quad [y] \mapsto [(y+\overline{y})/2].
$$

\begin{prop}\label{L is L}
$\alpha_*\co L_0(\ztt) \to
\mcL^{\ev}(\zt, 2)$ is an isomorphism.
\end{prop}
\begin{proof}
$ \ztt\xrightarrow{\alpha} \,_{2}(\Q [t]/\zt)$ and $
\ztt\xrightarrow{\beta}
\,_{2}(\Q[t]/2\Z[t]) =\Z_2[t]$ are isomorphisms.  Here $\,_{2}A =
\{a \in A : 2a = 0\}$ for an abelian group $A$.
\end{proof}

We  view $$\widehat H_n(\Z_2; R) = \frac{\{x\in R : 
\overline{x}=(-1)^{n+1}x \} }
{\{y-(-1)^n\overline{y} :\;\; y\in R\}}.
$$
     as a $R/2R$--module by the rule  $[a] \circ [x] := [ax\bar a]$.

We will need the following three classes of quadratic forms.
For each $p,g\in \widehat H_n(\Z_2;R)$ (or even $p,g\in R$ representing
classes in $\widehat H_n(\Z_2;R)$),  define
     $[P_{p, g}]\in L_{2n}(R)$  by the nonsingular
       $(-1)^n$--quadratic form:
\begin{equation*}
P_{p, g} = \left(R^2,
         \begin{pmatrix}
0&1\\
(-1)^n & 0
\end{pmatrix} ,
         \begin{pmatrix}
p\\
g
\end{pmatrix}  \right).
\end{equation*}
The notational conventions here are similar to those described after
Definition \ref{intronpg}.
The symmetric form $\lambda$ is determined by the matrix
$(\lambda(e_i,e_j))$  and
the  quadratic refinement $\mu$
    is specified by the vector $(\mu(e_i))$.

Next assume  $p, g\in R$
are in $\ker(1-*)$. Assume
$R$ is
torsion free as an abelian group. Define $[\mcP_{p, g}]\in \mcL^{\ev}(R,
2)$ and $[\mcN_{p, g}] \in \mcL(R, 2)$ by the quadratic linking forms:
\begin{align*}
\mcP_{p, g}  =&\left(( R/2R)^2,
\begin{pmatrix} 0 & 1/2 \\
1/2& 0
\end{pmatrix},
      \begin{pmatrix} p \\
g
\end{pmatrix}\right)\\
\mcN_{p, g} =&\left(( R/2R)^2,
\begin{pmatrix} p/2 & 1/2 \\
1/2& 0
\end{pmatrix},
\begin{pmatrix} p/2 \\
g
\end{pmatrix}\right)
\end{align*}

    Note that $ \mcN_{2p,
g} = \mcP_{p, g}$  and that  $[P_{p, g}]$ maps to $ [\mcP_{p, g}]$
under the composite
$$
L_2(R) \to L_2(R/2R) = L_0(R/2R) \xrightarrow{\alpha_*} \mcL^{\ev}(R,2).
$$

\begin{lemma}[Formal Properties of $P_{p,g}$, $\mcP_{p,g}$, and $
\mcN_{p,g}$]$\phantom{99}$
\label{formalproperties}
\rm
\begin{enumerate}

\item  $ P_{p,g} \cong P_{g,p}; \quad
    \mcP_{p,g} \cong
\mcP_{g,p}
$\quad  (
$\cong
$  means ``is isometric to").

\item $[P_{p_{1}, g}\oplus P_{p_{2}, g}] = [ P_{p_{1}+p_{2}, \,g}]; \quad
[\mcN_{p_{1}, g}\oplus \mcN_{p_{2}, g}] = [ \mcN_{p_{1}+p_{2},
\,g}]$.

\item $2[\mcN_{p, g}] = [\mcP_{p, g}]; \quad 2[\mcP_{p, g}] = 0;
\quad 2 [P_{p,g}]=0.$

\item \label{ps-p}$[P_{pg\overline{p}, g} ] = [P_{p,g}];\quad
[\mcP_{pg\overline{p}, g} ] = [\mcP_{p,g}]$.

\item\;
       $[\mcN _{\overline{x}px,\; g}] = [\mcN _{p, \; xg\overline{x}}] ;\quad [P
_{\overline{x}px,\; g}] = [P _{p, \; xg\overline{x}}]\quad $ for $x \in R$.

\item
       $[P_{p, g}]  = [P_{pg, 1}]$, if $R=\ztt,$ or if $R =
\Z[t]$ or $\Z[t^-]$ and $n$ is odd.

\item  $V_k([\mcN_{p,g}])= [\mcN_{p( t^k),
g(t^k)}]$ in $\mcL(R[t^\pm], 2).$

\item $[\mcN_{p, g+h}]-[\mcN_{p,g}]-[\mcN_{p,h}]
=[\mcP_{pg\overline{p},h}].$
\end{enumerate}
\end{lemma}

\begin{proof}

(1)\qua Replace the standard basis $\{e_1,e_2\}$ by $\{e_2,(-1)^n e_1\}$.

(2)\qua Let $\{e_1,e_2\}, \{e_3,  e_4\}$ be the standard bases for $P_{p_1, g}$
and $P_{p_2, g}$. A subLagrangian for $m=  P_{p_1, g}\oplus P_{p_2, g}$ is
$ S = \langle e_2+e_4\rangle.$  Then $\{[e_1+e_3], [e_2]\}$ is a
basis for $S^\perp/S$.
    The subLagrangian construction satisfies $m_S
=  P_{p_1+p_2, \,g}    $ proving the result.
The proof for $\mcN_{p_{1}, g}\oplus \mcN_{p_{2}, g}$ is the same.

(3)\qua is immediate from (2) and the definitions, since $2[p/2] = [p]$ in $\widehat H_1(\Z_2;R)$.

(4)\qua By (2) and (3) it suffices to note that $\langle e_1 + pe_2 \rangle$
is a  Lagrangian for $P_{pg\overline{p}+p, g}$.

(5)\qua A Lagrangian for $\mcN_{\overline{x}px, g}\oplus \mcN_{-p,
-xg\overline{x}}$ is
$\langle e_1+\overline{x}e_3, xe_2+e_4\rangle $. The proof of the second
part is similar.

(6)\qua We first consider the cases of $R = \Z[t]$ with $n$ odd and $R=
\ztt$. By (1) and (2), it is suffices
to consider the cases $p=t^k, g=t^l, $ for some
$k, l \geq 0$.
  From (5) above, $[P_{t^k, t^{l+2m}}]=[P_{t^{k+2m}, t^l}]$, so we may
assume $l$ is 0 or 1.  The case $l = 0$ is immediate, so we assume
$l = 1$ and proceed by induction on $k$.  If $k$ is even, the result
follows from (5) and (1).
If
$k = 2i-1$ is odd, then
$[P_{t^{2i-1},t}] = [P_{t^{i-1},t}] = [P_{t^i,1}] = [P_{t^{2i},1}]$;
the first equality comes from replacing  the
standard  basis $\{e_1, e_2\}$   by
$\{e_1+t^{i-1}e_2, \;e_2\}$, the second equality from the induction
hypothesis, and the third from (4).

The proof in the case $R = \Z[t^-]$ and $n$ odd is similar; it
involves replacing $t$ by $t^2$ everywhere in the above paragraph.

(7)\qua This is immediate from the definitions.

(8)\qua  Let $m= \mcN_{p,g+h}\oplus\mcN_{-p,g}\oplus\mcN_{
-p,h}$. A subLagrangian for $m$ is
$S=\langle e_1+e_3, e_2+e_4+e_6\rangle$.
An arbitrary element $\sum_{i=1}^6 a_ ie_i \in (R/2R)^6$  is in $S^\perp$
if and only if: $a_1p +a_2 +a_3p+a_4=0=a_1+a_3+a_5$. Therefore, $x\in
S^\perp $ if and only if $x=a_1(e_1+e_3)+a_2(e_2+e_4)+a_5(e_3+pe_4+e_5)
+a_6 e_6$. Hence $S^\perp/S = \langle [e_3+pe_4+e_5], [e_6]\rangle$.  Thus
   $[m] = [m_S]=[\mcP_{pg\bar{p},
h}]$.
\end{proof}

\begin{definition}  \label{tate}    Let $R$ be a ring with involution.  Define  maps
\begin{enumerate}
\item  
$P_2 = P_2^R \co  \widehat H_1(\Z_2; R) \to L_2(R), \qquad  [p] \mapsto [P_{p,1}],$

\item  
$Q_0 = Q_0^{R[t^{-}]} \co  t^2 \widehat  H_1(\Z_2; R[t^{-}]) \to NL^-_0(R), \qquad  [t^2 p] \mapsto [P_{tp,t}].$
  
\end{enumerate}

  Note that 
$
\varepsilon_{0*}[P_{tp,t}] = [P_{0,0}] = 0$, so $Q_0$ takes its values in $N^-L_0(R) \subset L_0(R[t^{-}])$.   Both $P_2$ and $Q_0$ are homomorphisms by Lemma \ref{formalproperties}(2).
\end{definition}

\begin{lemma} \label{trivial}
 Let $R$ be a commutative ring with trivial involution.  
\begin{enumerate}
\item  If $p \in R[t^\pm ]$ satisfies $\overline p = p$, then $P_2^{R[t^\pm ]}([p^2] - [p]) = 0$.
\item  If $p \in R[t^- ]$ satisfies $\overline p = p$, then $Q_0^{R[t^- ]}([(t^2p)^2] - [t^2p]) = 0$.
\item  $P_2^{R[t]}$ is a map of $\mcV$--modules.  $P_2^{R[t^-]}$ and $Q_0^{R[t^-]}$ are maps of $\mcV_{\od}$--modules.  (The Verschiebung algebras $\mcV$ and $\mcV_{\od}$ were defined in the introduction.)
\end{enumerate}
\end{lemma}

\proof  Parts (1) and (2) follow from Lemma \ref{formalproperties}(4).  The assertions concerning $P_2$ in Part (3) are clear.  Finally, note
\begin{align*}
V_{2k+1}Q_0[t^2p(t^2)] & = 
[P_{t^{2k+1}p(t^{4k+2}),t^{2k+1}}]  \\
& = 
[P_{t^{4k+1}p(t^{4k+2}),t}] \qquad  \text{by Lemma \ref{formalproperties}(5)}\\
& = Q_0 V_{2k+1}[t^2p(t^2)].\tag*{\qed}
\end{align*}

We now specialize to  $R = \Z$ and $\Z_2$.  We will abuse notation (somewhat) in the statement of the theorem below by maintaining the names $P_2$ and $Q_0$ for a factorization through a quotient.

\begin{thm}\label{3n3}  
\begin{enumerate}

\item $\eps_{0 *}\co L_0(\Z[t]) \xrightarrow{\cong} L_0(\Z).$

\item   $ P^{R[t]}_2\co \widehat H_1(\Z_2;R[t])/(V_2-1) \widehat H_1(\Z_2;R[t]) \xrightarrow{\cong}
L_2(R[t]), $ if
$R = \Z$  or $\Z_2$.

\item $ P^{\Z[t^-]}_2\co \widehat H_1(\Z_2;\Z[t^-])/\{p^2 - p : p \in \widehat H_1(\Z_2;\Z[t^-])\} \xrightarrow{\cong}
L_2(\Z[t^-])$.

\item  $ Q^{\Z[t^-]}_0\co  t^2 \widehat H_1(\Z_2;\Z[t^-])/\{p^2 - p : p \in t^2 \widehat H_1(\Z_2;\Z[t^-])\} \xrightarrow{\cong}
N^-L_2(\Z[t^-])$.

\end{enumerate}

The isomorphism in (1) is of abelian groups, in (2) of $\mcV$--modules, and in (3) and (4) of $\mcV_{\od}$--modules.
\end{thm}

To prove this, we will need the following lemma, which is similar to 
\cite[Proposition 41.3(v)]{RHK} and \cite[Proposition 2.11(ii)]{CR}.

\begin{lemma}  \label{3n4} Let $R$ be  a  principal ideal domain with
involution.  Any $\xi \in \ker( \eps_{0*}\co L_{2n}(R[t^{\pm}])\to L_{2n}(R)) $ can be
represented by a
$(-1)^n$--quadratic  form $(M[t],\lambda, \mu)$ for which there is a free
$R $--summand $L\subset M$, such that $L[t]= L[t]^\perp$ (a ``symmetric
Lagrangian").
\end{lemma}

\begin{proof} \rm By Higman Linearization (\cite[Lemma 3.6ab]{CK}, also \cite[Proposition 5.1.3]{RES}), 
extended to the case $t \mapsto -t$,  
one can represent
$\xi$ by a form $(M[t],\lambda, \mu)$ $= (M[t], \lambda_0+t\lambda_1, \mu_0+
t\mu_1).$  Here $M$ is a finitely generated
free
$R$--module; $M[t]=R[t]\otimes_R M$; and $\lambda_0(x,y), \lambda_1(x,y)$ are $R$--valued.     
A similar interpretation holds  for $\mu_0 +t\mu_1$.

It follows that $(M,\lambda_0,\mu_0)$ is a $(-1)^n$--quadratic form and $(M,\lambda_1,\mu_1)$
is an $\eps (-1)^n$--quadratic form where $t \mapsto \eps t$. 
Since
$\lambda_0$ and
$\lambda_0+t\lambda_1$ are nonsingular,  the $R$--map
$\nu = (\ad \lambda_0)^{-1}\circ \ad \lambda_1\co   M \to M$ is nilpotent and satisfies:
\[
\lambda_1(x,y)=\lambda_0(\nu x,y)=\eps \lambda_0(x,\nu y), \quad
\text{for all }  x,y \in
M.
\]

Choose $k\geq 0$ so that $\nu^{k+1} = 0, $ but $\nu^k\neq 0$. First assume $\nu >0$.  Let $$\overline{\nu^k(M)} = \{x \in M : ax \in \nu^k(M) \text{ for some } a \in R - 0\}.$$   
This is a summand of $M$.
 Pick a basis element $e_1\in \overline{\nu^k M}$. Then $\nu(e_1)=0,$ and $\lambda_0(e_1,e_1)=0.$

Set $m= (M,\lambda_0),\; S_1= Re_1$.  Then
$S_1$ is a subLagrangian for $m$ and $\nu$ induces a
nilpotent map
$\nu_{S_1}\co  S_1^\perp/S_1\to S_1^\perp/S_1$.
Repeat this step on
the  subLagrangian construction $m_{S_1} =(S_1^\perp/S_1,
(\lambda_0)_{S_1})$, getting $e_2$, etc., until one obtains a basis
$\{e_1, e_2,\dots,e_m\}$ for a summand $S$ of $M$ satisfying:
\[
S\subset S^\perp, \quad \nu(S^\perp)\subset S.
\]
The Witt class of $m_S  =(S^\perp/S, (\lambda_0)_S)$ is zero
because
$\xi\in \ker(\eps_{0*})$. So by  adding a hyperbolic form,
if necessary, to
the original $(M[t], \lambda, \mu)$, we may as well assume $m_S$ is
hyperbolic. We can therefore find additional elements,
      $e_{m+1},\dots, e_k\in
S^\perp$, whose images in $S^\perp/S$ form a basis for a Lagrangian of
$m_S$. It follows that
$L=< e_1,\dots , e_k>$ is a summand of $M$ satisfying:
\[
L^\perp = L \;(\text{relative to }\lambda_0); \quad \nu(L)\subset L.
\]
Therefore $L[t] = L[t]^\perp $ (relative to $\lambda)$.
\end{proof}

\begin{proof}[Proof of Theorem \ref{3n3}]

Let $R= \Z$ or $ \Z_2$.  Let $n = 0$ or $1$.  Let $\xi \in N^{\pm}L_{2n}(R) = 
\ker(\eps_{0*} \co L_{2n}(R[t^\pm] \to L_{2n}(R))$. 
By Lemma \ref{3n4}, we may write
$\xi = [M[t], \lambda, \mu]$, where $M$ is a free $R$--module containing a free summand $L$ for which
$L[t]=L[t]^\perp$. Since $\lambda$ is an even form, any $R$--basis $\{e_1,\dots , e_k\}$ for $L$ extends to an
$R[t]$--symplectic basis $\{e_1,\dots, e_k, f_1,\dots , f_k \}$ for $(M[t],\lambda)$.
So
\[
\xi=\sum_{i=1}^k~[P_{p_i, q_i}]\; \text{  where } p_i =\mu(e_i),\;
q_i=\mu(f_i).
\]
We use this expression for $\xi$ repeatedly below. 

Our first step is to show that the maps $P_2^{R[t^{\pm}]}$ and $Q_0^{\Z[t^-]}$ of Theorem \ref{3n3} are surjective.  We start with $P_2^{R[t^{\pm}]}$.  Note that $P_2^{R[t^{\pm}]}[1] = [P_{1,1}]$ generates the image of $L_2(R)$, so we just need to show that every $\xi \in N^\pm L_2(R)$ is in the image of $P_2$.  By the first paragraph of this proof and parts (6) and (2) of Lemma \ref{formalproperties}, 
$$
\xi = \sum_{i=1}^k~[P_{p_i, q_i}] = \sum_{i=1}^k~[P_{p_i q_i,1}] = P_2^{R[t^{\pm}]}( \sum_{i=1}^k p_i q_i).
$$

To see $Q_0^{\Z[t^-]}$ is surjective, consider $\xi = \sum_{i=1}^k~[P_{p_i, q_i}] \in N^-L_0(\Z)$ as before.  Note $p_i, q_i \in \widehat H_0(\Z_2; \Z[t^-]) = t \Z_2[t^2]$.  Therefore, write $p_i \equiv r_i t \overline r_i \pmod {2\Z[t^-]}$ and $q_i \equiv s_i t \overline s_i \pmod {2\Z[t^-]}$ where $r_i,s_i \in \Z[t^-]$.  Then, with the help of Lemma \ref{formalproperties}(5), one sees
$$
\xi = \sum_{i=1}^k~[P_{p_i, q_i}] = \sum_{i=1}^k~[P_{r_i t \overline r_i, s_i t \overline s_i}] =
 \sum_{i=1}^k~[P_{\overline s_i r_i t \overline r_i s_i, t }] = Q_0^{\Z[t^-]}(  \sum_{i=1}^k \overline s_i r_i t^2 \overline r_i s_i).
$$

Now we need to prove injectivity of the four maps.  We start with $\varepsilon_{0*}$.  By the first paragraph of the proof, every $\xi \in \ker(\varepsilon_{0*} \co L_0(\Z[t] \to L_0(\Z))$ can be represented by a quadratic form $(M[t],\lambda, \mu)$ where $M[t]$ has a free summand $L[t]$ so that for all $x \in L[t]$, one has $\lambda(x,x) = 0$.  But $2\mu(x) = \lambda(x,x) \in \Z[t]$, so $\mu|_{L[t]} = 0$.  Thus $L[t]$ is a Lagrangian and so $\xi = 0$.  Hence $\varepsilon_{0*}$ is injective.

For a field $F$ of characteristic 2, there is a well-defined epimorphism called the Arf invariant (see \cite{Arf}, \cite{Sc}) 
$$
\Arf\co  L_{2n}(F) \to F/\{x^2 + x : x \in F\}, \quad [M,\lambda,\mu] \mapsto \sum_{i = 1}^{k} \mu(e_i)\mu(f_i),
$$
where $\{e_1, \dots e_k, f_1, \dots f_k\}$ is a symplectic basis for $(M,\lambda)$.  Applying this to the function field $\Z_2(t)$ and then restricting to the subring $\Z_2[t]$, one obtains an epimorphism
$$
\arf \co L_{2n}(\Z_2[t]) \to \Z_2[t]/\{x^2 + x : x \in \Z_2[t]\}.
$$
(Note, there is a subtle point here -- one needs that $\Z_2[t]/\{x^2 + x : x \in \Z_2[t]\} \to \Z_2(t)/\{x^2 + x : x \in \Z_2(t)\}$ is injective.  This can be shown using the fact that $\Z_2[t]$ is a PID and hence integrally closed.)

We now return to the proof of injectivity of $P_2^{R[t]}$ where $R = \Z$ or $\Z_2$.  Let $r : L_2(R[t]) \to L_2(\Z_2[t])$ be mod 2 reduction or the identity map.  Then $P_2^{R[t]}$ is injective since
$
\arf \circ ~ r \circ P_2^{R[t]}
$
is the identity map.

To show that $P_2^{\Z[t^-]}$ is injective we consider mod 2 reduction $r\co L_2(\Z[t^-]) \to L_2(\Z_2[t])$.  We will show
$$
\arf \circ ~r \circ P_2^{\Z[t^-]}  :  \Z_2[t^2]/\{x^2 - x : x \in \Z_2[t^2]\} \to  \Z_2[t^2]/\{x^2 - x : x \in \Z_2[t]\}
$$
is injective.  Suppose $(\arf \circ ~r \circ P_2^{\Z[t^-]})[p] = [p] = 0$ for some $p \in \Z_2[t^2]$.  This last equality means $ p = x^2 - x$ for some $x \in \Z_2[t^2]$.  Since $x^2$ and $x^2 - x$ are both in $\Z_2[t^2]$, we see $x \in \Z_2[t^2]$ and thus $[p] = 0 \in \Z_2[t^2]/\{x^2 - x : x \in \Z_2[t^2]\}$.  Hence $P_2^{\Z[t^-]}$ is injective.

To show that $Q_0^{\Z[t^-]}$ is injective we consider mod 2 reduction $r\co L_0(\Z[t^-]) \to L_0(\Z_2[t])$.  We will show
$$
\arf \circ ~r \circ Q_0^{\Z[t^-]}  : t^2 \Z_2[t^2]/\{x^2 - x : x \in t^2\Z_2[t^2]\} \to  \Z_2[t^2]/\{x^2 - x : x \in \Z_2[t]\}
$$
is injective.  Suppose $(\arf \circ ~r \circ Q_0^{\Z[t^-]})[t^2p] = [t^2p] = 0$ for some $p \in \Z_2[t^2]$.  This last equality means $ t^2 p = x^2 - x$ for some $x \in \Z_2[t^2]$.  Since $x^2$ and $x^2 - x$ are both in $\Z_2[t^2]$, we see $x \in \Z_2[t^2]$.  Set $x = c+y$ where $c \in \Z_2$ and $y \in t^2 \Z_2[t^2]$.  Clearly $y^2 - y = (x^2 - x) + (c^2 - c) = x^2 - x = t^2 p$.  Thus $[t^2p]=0$ in the domain.  
 Hence $Q_0^{\Z[t^-]}$ is injective.
\end{proof}

We can now calculate $\mcL^{\ev}(\Z[t^\pm],  2)$.

\begin{thm}\label{even}
\begin{enumerate}
\item  There
is an isomorphism of $\mcV$--modules,
\begin{multline*}
\mcV/\langle 2, V_2-1\rangle \cong\frac{t\Z_2[t]}{( V_2-1)t\Z_2[t] }
\overset{\mcP}{\cong}
       \widetilde{\mcL}^{\ev}(\Z[t], 2); \\ \mcP[p]  =\left[\ztt^2,
\begin{pmatrix} 0 & 1/2 \\
1/2& 0
\end{pmatrix}
\begin{pmatrix} p \\
1
\end{pmatrix}\right].
\end{multline*}
\item   \quad $\widetilde{\mcL}^{\ev}(\Z[t^-], 2)=0.$

\end{enumerate}
\end{thm}
\begin{proof}
(1)\qua By Proposition \ref{L is L},
$
\alpha_*\co  \widetilde L_0(\Z_2[t])\to  \widetilde \mcL^{\ev}(\Z[t], 2)
$
       is an
isomorphism.  Then $(\alpha_*)^{-1} \circ \cal P$ is the isomorphism of
Theorem \ref{3n3}(2).

(2)\qua Let $m=(M, b, q)$ be a \qlf of
exponent
$2$ with even type over $\Z[t^-]$. Then $M$ admits a symplectic basis over
$\ztt$. 
This means $m$ is an orthogonal sum of  terms of
the form $\mcP_{p,g}$ , where $p,g\in\Z[t^2]$.  By parts (2) and (4) of Lemma
\ref{formalproperties},
$[m]$ is a sum of terms of the form $[\mcP_{t^{4n+2},1}]$.

But $[\mcP_{t^{4n+2},1}]=0$ in $\widetilde{\mcL}(\Z[t^-],2)$, because
$\langle e+t^{2n+1}f\rangle$ is a Lagrangian for
$\mcP_{t^{4n+2},1}$. To see
this,
       we compute:
\[
q(e_1+t^{2n+1}e_2)=
t^{4n+2}+t^{4n+2}+\{1+*\}(t^{2n+1}/2)=\{1+*\}(t^{2n+1}/2)=0.
\]
This proves $\widetilde{\mcL}^{\ev}(\Z[t^-], 2)=0.$
\end{proof}

\section{Characteristic elements and the proof of Theorem
\ref{devissage}}\label{S3}

In this section we    write $R$ for $\Z[t]$, $\Z[t^-]$, or $\Z$.
For any $R$--module $M$ and any $k\in \N $  we write 
$$
{}_k M = \{x\in M :kx=0\}.
$$
    Let $N  $ be a submodule of a length one $R$--module $M$ of  exponent
$2^n$, for some $n$.  The \emph{closure} of $N$ in $M$ is:
\[
\overline{N} = \{x\in M:px\in N \text{  for some } p\in R\smallsetminus
2R\}.
\]

By Corollary
\ref{colengthcor} of the Appendix, $\overline{N}$ is the unique  smallest
colength one submodule of
$M$ which contains $N$.  We will use this concept below.

The goal of this section is to prove:
\begin{thm}\label{ses2q}
Let $R=\Z[t^\pm]$  or $\Z$.  There is a short exact sequence:
\[
0\to \mcL(R, 2 )\to
\mcL(R, \langle 2 \rangle)\xrightarrow{Q}\Z_2 \to 0
\]
           Moreover,  $ \mcL(R,
2^n)\to\mcL(R, 2^{n+1}) $ is an
isomorphism if
$n\geq 2$. 
\end{thm}
We will then prove Theorem \ref{devissage} as a corollary of
Theorem \ref{ses2q}.

\begin{cce}

Fix an integer $n\geq 1$ and a \qlf $m=(M,b,q) $ for which $2^nM=0$.
We  construct elements
$v^{(n)}_0(m), v^{(n)}_1(m)\in{}_{2}M$. They depend on $M, b,$ and
$n$.

Recall that  $\widehat H_1(\Z_2; R)$ is an $R/2R$--module via $[a] \circ [x] = [ax\bar a]$.
Consider  the $R$--linear map
\[
          \phi\co  M\to \widehat H_1(\Z_2; R),\quad     \phi(x) =
\{1+*\}2^{n-1}b(x,x)=2^nq(x)\in \widehat H_1(\Z_2;R).
\]
Here $ \{1+*\}\co \frac{1}{2}R/R\to \widehat H_1(\Z_2; R)$ is given by
 $[x] \mapsto [x+\bar x] \in \widehat  H_1(\Z_2; R)$. It is bijective if
$R=\Z[t]$ or $\Z$.  

The map $\phi$ measures  the failure of $m$ to be of even type.

As an  $R/2R$--module, $\widehat H_1(\Z_2; R)$ is free on the basis
$\{[1], [t]\}$ if
$R=\Z[t]$ and on the basis $\{[1]\}$ if $R=\Z[t^-]$ or $\Z$.
Therefore, if
$R=\Z[t]$, there are
$R$--maps $\phi_0, \phi_1\co  M\to R/2R$, uniquely
specified by the rule
\[
\phi(x) = \phi_0(x)\circ[1]+\phi_1(x)\circ [t].
\]
If $R=\Z[t^-]$ or $\Z$, we get a single map $\phi_0 \co  M\to R/2R$, uniquely specified by
\[
\phi(x) = \phi_0(x)\circ[1] .
\]
We define $\phi_1$ to be $0$ in these cases.

The map $ \ad b\co M \to M^\wedge$ restricts to an isomorphism
\[ \;{}_{2}M \xrightarrow{\cong}
{}_{2}(M^\wedge) =\Hom_R (M,
\frac{1}{2}R/R).
\]  
We will change the target using the isomorphism $\{2\}\co \frac{1}{2}R/R \xrightarrow{\cong} R/2R $ sending $x+R$ to
$2x+ 2R$.

Therefore there are elements $ v_i = v^{(n)}_i(m)\in {}_{2}M$, defined by
$$
v_i = (\ad b)^{-1}(\{2\}^{-1} \circ \phi_i).
$$
In other words
\[
\{2\}b(v_i,x) =\phi_i(x)  \text{ for all $x \in M$. }
\]

 If $R=\Z[t^-]$ or $\Z$, then  $v_1=0$.

       We conclude that $v_0, v_1 \in {}_2M$ are characterized by the 
fact that for all $x\in M$:
\begin{align}
\label{veqn1}
           \{2\}2^{n-1}b(x,x)  &= (\{2\}b(v_0, x))^2 + (\{2\}b(v_1,
x))^2t \in \ztt,
\text{ if } R=\Z[t],\\
\label{veqn2}
\{2\}2^{n-1}b(x,x)  &= (\{2\}b(v_0, x))^2\in \Z_2[t^2], \text{ if }
R= \Z[t^-],\\
\label{veqn3}
\{2\}2^{n-1}b(x,x)  &= (\{2\}b(v_0, x))^2 = \{2\}b(v_0, x)\in \Z_2, \text{
if }  R=\Z.
\end{align}
\end{cce}

\begin{prop}\label{onvi}
The  construction above satisfies the following properties:
\begin{enumerate}
\item If $m= m'\oplus m'' $\;(orthogonal direct sum), then
\[v_i^{(n)}( m'\oplus m'')
= v_i^{(n)}(m')\oplus v_i^{(n)}(m'').
\]
\item  Suppose $L\subset M $ and $2^{n-1}b(x,x)=0$,   for all $x\in L$.
Then
$v_i^{(n)}(m)\in L^\perp $.  So if $L$  is a
Lagrangian for $m$, then $v_i^{(n)}(m)\in L$, and $q(v_i^{(n)}(m))=0$ for
each i.
\item   \quad $v_i^{(n)}(m)\in \overline{2^{n-1}M}$, for $ i=0,1$.
\end{enumerate}

\end{prop}
\begin{proof}

(1)\qua This is obvious from the definitions.

(2)\qua This is clear from \eqref{veqn1}, \eqref{veqn2}, and
\eqref{veqn3}, because
$\phi(L) = 0$.

(3)\qua   In general, $(2^kM)^\perp = {}_{2^k}M$. This implies
$({}_{2^k}M)^\perp = (2^k M)^{\perp\perp}= \overline{(2^kM)}$, by
Proposition \ref{colengtheq} in the Appendix. So we must only
show that
$v_i\in({}_{2^{n-1}}M)^\perp$, for
$i=0,1$. But this is clear from
          (2) above, if we set $L={}_{2^{n-1}}M$.
\end{proof}

There is one crucial case in which we can slightly strengthen
Proposition \ref{onvi}(3), which we describe now.
A \qlf $(M, b, q)$  is {\em  irreducible } if for all $x\neq 0$ in
$M$, we have $q(x)\neq 0$.
\begin{lemma}\label{vin2}  Let $m=(M, b, q)$ be an
irreducible \qlf over $R$, with $ 4M=0$.  Then $v_i^{(2)}(m)\in
2M$, for each $i$.
\end{lemma}

\begin{proof}
We assume $R=\Z[t^{\pm}]$; in the $R = \Z$ case there is nothing to prove by 
 Proposition \ref{onvi}(3). For any
$p\in \zt  \smallsetminus 2\zt$,  let
$M(p)$ denote the
$\zt$--module generated by two elements  $\phi,\tau$, subject only to the two
relations:
\[2\tau=0,\qquad2\phi = p\tau.
\]
Note that $M(p)$ can also be viewed as the submodule $\langle [p], [2]\rangle$ of $\Z_4[t]$.
Hence by 
Corollary \ref{lengthcor} $M(p)$ has length one.

Its dual,
$M(p)^\wedge = \Hom_R(M,\Q R/R)$, is generated by  two homomorphisms $\Phi,  T  $ defined
by:
\[
\Phi(\phi)=\frac{p}{4},\quad  \Phi(\tau)=\frac{1}{2}= T(\phi),\quad
T(\tau)=0.
\]
Note $M(p)\cong M(p)^\wedge$.

Since $4M = 0$, multiplication by 2 gives a monomorphism 
$\times 2\co  M/{}_{2}M
\to {}_{2}M$  of
finitely generated torsion-free modules over the principal ideal domain $\ztt$.
By the structure theorem for such modules, we can find
bases $\{\tau_1,\tau_2,\dots, \tau_l\}$ for ${}_{2}M$, and $ \{[\phi_1],
[\phi_2],\dots, [\phi_k]\} $ for $M/{}_{2}M$ over $\ztt$, (where
$\phi_i\in M, i=1,\dots, k$), so that $2\phi_i=p_i\tau_i$ for some $p_i\in
\zt \smallsetminus 2\zt$, whose class  $[p_i]\in \ztt$  divides $[p_{i+1}]$,  for
$0<i<k$. Note that $0\leq k\leq l.$

It follows that $ M = M(p_1)
\oplus M(p_2)\oplus\dots  M(p_k)\oplus N$, where $N\cong \Z_2[t]^{l-k}$ is the submodule
with basis $\{\tau_{k+1}, \dots, \tau_l\}$. Therefore, $M^\wedge $ is spanned by
$\Phi_1,\dots, \Phi_k, T_1,\dots, T_l$, and $N^\wedge $ is spanned by
$T_{k+1},\dots ,T_l$.

Now rank$_{\ztt}(\widehat H_1(\Z_2; R))= 1$ or $2.$ Since $m$ is
irreducible,
\[
4q\co  M/{}_{2}M\to \widehat H_1(\Z_2;R)\]
is an
$\ztt$--linear monomorphism. So
$k\leq 2$ if $R=\Z[t]$; $k\leq 1$ if $R=\Z[t^-]$.
The rest of the proof breaks into cases depending on the value of k.

\medskip

\bf Case One\rm\qua Assume $ k=0$.  By Proposition \ref{onvi}(3), $v_i^{(2)}(m) \in \overline{2M}$, but $\overline{2M} = 2M$ since $2M = 0$.  

\medskip

\bf Case Two\rm\qua Assume $k=1$.   We
have
$M=M(p_1)\oplus N$. Because
$b$ is nonsingular, $(\ad b)(\phi_1)\equiv\Phi_1 \pmod {{}_{2}(M^\wedge)}$, 
noting that $M/{}_2M$ and $M^\wedge/{}_2(M^\wedge)$ are both free and rank 1 over $\Z_2[t]$ with bases $\phi_1$ and $\Phi_1$.
Therefore, $2b(\phi_1, \phi_1)=\frac{p_1}{2}$.
We have:
\[
p_1=\{2\}2b(\phi_1, \phi_1) =4q(\phi_1)= q(p_1\tau_1) =
(p_1)^2q(\tau_1)\in \ztt.
\]
So $p_1q(\tau_1) =1$ in $\ztt$. In particular, $[p_1]=1\in \Z_2[t]$, and so
$2M=\overline{2M}$ since $M = \Z_4[t] \oplus \Z_2[t]^j$. So again Proposition \ref{onvi}(3)  proves the result.

\medskip

\bf Case Three\rm\qua Assume $k=2$.
Hence $R=\Z[t]$ and  $ M=
M(p_1)\oplus M(p_2)\oplus N.$

\begin{claim}
 $[p_1]=[p_2]\in \ztt$.
\end{claim}

\begin{proof}[Proof of Claim] Write $[p_2]=[p_1r]$ where  $ r\in \zt$. Then
let  $A\co \overline{2M}\cong\overline{2M^\wedge}$  be the  restriction of
the isomorphism $\ad b\co M\to M^\wedge$. Let
$\left(\begin{smallmatrix}
\alpha &\beta \\
\gamma &\delta 
\end{smallmatrix}
\right) \in GL_2(\ztt)
$
be the matrix of $A$ relative to the  bases $\{\tau_1, \tau_2\} $ and
$\{T_1,  T_2\}$ for $ \overline{2M}$ and $\overline{2M^\wedge}$
respectively. We then have
\begin{multline*}
\frac{p_1}{2}\gamma  = p_1\gamma T_2(\phi_2) =
p_1A(\tau_1)(\phi_2)=p_1b (\tau_1,\phi_2) = p_1b(\phi_2, \tau_1) =\\
b(2\phi_2, \phi_1)=p_2A(\tau_2)(\phi_1) =
p_1r\beta T_1(\phi_1)=\frac{p_1}{2}r\beta \in \frac{1}{2}\zt/\zt.
\end{multline*}
Therefore $\gamma =r\beta $.  The matrix
$\left(\begin{smallmatrix}
\alpha &\beta \\
r\beta &\delta 
\end{smallmatrix}
\right)
$
is therefore nonsingular, $\alpha \delta =1+r\beta ^2$, and  consequently $\delta $ is
relatively prime to $[r]$.

On the other hand, 
\[
p_1r\delta = p_2\delta = p_2\{2\}\frac{\delta }{2} = \{2\}b(p_2\tau_2,\phi_2)
        = 4q(\phi_2) = q(p_2\tau_2) = p_1^2r^2q(\tau_2) \in \ztt.
\]
Therefore $\delta =p_1r \,q(\tau_2)$, and so $[r]$ obviously divides $\delta $.
This
implies $[r]=1\in \ztt$, proving the claim.
\end{proof}

Returning to the proof  of Case Three, we write $M=M(p)\oplus M(p)\oplus
N$. Note that $4q(\phi_i) =p^2q(\tau_i), i=1,2$. Therefore
$4q(M)\subset  p\circ\widehat H_1(\Z_2;\zt),$ and so $\{2\}b(v_i,
M)\subset p\ztt$,
which tells us  $ v_i\in pM$. But, by construction of $M$ we see
$pM\cap\overline{2M}= 2M$, so by Proposition
\ref{onvi}(3) again, $v_i\in 2M,$ for $i=0, 1$. This completes the proof.
\end{proof}

\begin{remark} A  study of the above proof shows that $M(p)$
has the following two {\em cool} properties.  Each length one $\Z[t]$--module
of exponent four is isomorphic to a direct sum $M(p_1) \oplus \dots
M(p_k) \oplus \Z_2[t]^j$. If $[p]\not  = 0,1$,
then  $M(p)$ admits a nonsingular bilinear linking form, but no
quadratic linking form.
\end{remark}

Our proof of Theorem \ref{ses2q} is largely accomplished by the following
definition and the next two lemmas.

\begin{definition}(of $\qnz, \qnone:\mcL(R,2^n)\to \widehat
H_1(\Z_2;R)$)

Let $ m=(M, b, q)$ be a quadratic linking form over $R$ with $2^nM=0$. By Proposition
\ref{onvi}(3), the characteristic elements satisfy $v_i= v^{(n)}_i(m)\in
\overline{2^{n-1}M}$. Therefore if $n\neq 1$,
$b(v_i, v_i) = 0$ for each  $i$, so that $q(v_i)\in \widehat
H_1(\Z_2;R) $. (Recall
$\widehat H_1(\Z_2;R)$ is $\Z_2[t]$ or $ \Z_2[t^2]$ or $\Z_2$.)
We set
\begin{align*}
       Q^{(n)}_i(m) &= q(v_i)\in \widehat H_1(\Z_2;R), \text{ if } n > 1,\\
Q^{(1)}_i(m) &= q(v_i)\in \frac{1}{2}\{x \in R : x = \overline x\}/(1+*)R, \text{ if } n = 1.
\end{align*}
Then $Q^{(n)}_i$ is additive, and $Q^{(n)}_i(m) =0$ if $m$ admits a
Lagrangian, by Proposition \ref{onvi}(2). So each $Q^{(n)}_i$  is a 
well-defined homomorphism on $\mcL(R,2^n)$.

Obviously, if $R=\Z$ or $\Z[t^-]$,
$Q^{(n)}_1=0 $ for all $n$.
We define
\begin{equation}\label{qnz}
\mcL^{\eps}(R, 2^n)= \ker{\qnz}\cap\ker{\qnone}\subset \mcL(R, 2^n).
\end{equation}
\end{definition}

If $m$  is a \qlf of exponent $2^n$ with
even type, $\vnz(m)=\vnone(m)=0$.  Hence
\begin{equation}\label{linl}
\iota_n(\mcL^{\ev}(R, 2^n))\subset \mcL^{\eps}(R, 2^n).
\end{equation}

       \begin{lemma}\label{imq}
\begin{enumerate}
\item For all
$n\geq 3,\quad \mcL^{\eps}(R, 2^n)=\mcL(R, 2^n)$ .
\item Furthermore $Q^{(2)}_1=0$, and there is a short
exact sequence:
$$
0\to \mcL^{\eps}(R, 4)\to\mcL(R,
4)\xrightarrow{Q^{(2)}_0} \Z_2\to 0.
$$

\item Finally, for each $i$, $Q^{(1)}_i(\widetilde{\mcL}(R, 2))\subset
\widehat H_1(\Z_2;R)$.
\end{enumerate}
\end{lemma}
\begin{proof}
(1)\qua First suppose $n\geq 3$, and $[m]=[M, b, q]\in
\mcL(R, 2^n)$.  We  show  for each
$i$, that $ q(v^{(n)}_i(m)) =0$. This will prove
$\mcL^{\eps}(R, 2^n)=\mcL(R, 2^n)$.

By Proposition \ref{onvi}(3), $
v^{(n)}_i(m)\in\overline{ 2^{n-1}M}$. Therefore $pv^{(n)}_i(m)\in
2^{n-1}M$,  for some $p\in R$ such that $[p]\neq 0$ in $R/2R.$ Therefore
$p^2q(v^{(n)}_i(m)) \in q(2^{n-1}M)= 2^{2n-2}q(M) = 0$, because
       $2n-2\geq n+1$.  This   implies
$q(v^{(n)}_i(m))=0$.

(2)\qua  Next we prove $Q^{(2)}_1=0$ and $\im{(Q_0^{(2)})}\subset \Z_2$.

Let
$m=(M, b, q)$ be a
\qlf with $4M = 0$.  Let $v_i = v_i^{(2)}(m)$. We must
show that
$q(v_1)=0,$ and
$q(v_0)=0$ or $1.$   If we
choose, we may replace
$m$ by
$m_S =(S^\perp /S, b_S, q_S)$, where $S$ is any subLagrangian of $m$.
Moreover, $S=\overline{\langle x\rangle}$ \emph{ is} a subLagrangian, for
any
$x\in M$ for which  $q(x)=0$, by Corollary \ref{colengthcor} of the
Appendix. This means we may as well assume that $m$ is irreducible.
Therefore, by Lemma \ref{vin2}, we may write $v_0 = 2y_0, \;v_1=2y_1$, for
some $y_0, y_1$.  Define $b_{ij} =\{2\}b(v_i,y_j)= b_{ji} \in R/2R.$ 
Note 
$$
b_{ii} = \{2\}2b(y_i,y_i) = 4q(y_i) = q(v_i).
$$
From
\eqref{veqn1} we see
\[
b_{ii}=\{2\}2b(y_i,y_i)= b_{0i}^2 + b_{1i}^2t.
\]
(When $R=\Z[t^-]$ or $\Z$, these equations simplify to: $b_{00} = b_{0 0}^2$,
which immediately gives $Q^{(2)}_0(m) = 0 $ or $ 1, Q^{(2)}_1=0$).

But, in the ring $R/2R=\Z_2[t]$, the equations
\begin{align*}
b_{00}&= b_{00}^2 + b_{10}^2 t\\
b_{11}&= b_{10}^2 + b_{11}^2 t
\end{align*}
have only two solutions:  $ 0=b_{11}=b_{10};\quad b_{00}= 0$ or
$1$. For no polynomial of the form $p^2-p\in \Z_2[t]$ is  an odd degree
polynomial if
$\text{deg}(p)>0$. So $Q^{(2)}_0(m)\in
\{0, 1\}$,
$Q^{(2)}_1(m)=0,$ and $\mcL^{\eps}(R, 4) = \ker(Q^{(2)}_0)\cap
\ker(Q^{(2)}_1) =
\ker(Q^{(2)}_0)$.

To finish the proof of (2), we note that
$Q^{(2)}_0([m])=1
$ when
\[
M=\Z_4[t^\pm], \quad b(x, y) =
\frac{x\bar{y}}{4} , \quad q(x)= \frac{x\bar{x}}{4} \quad \text{ for all }
x,y\in M.
\]
          In this
case  $ v_0 = 2$.
Consequently, $Q^{(2)}_0(m) = q(2)= 1$. This proves (2).

(3)\qua This is
similar to the argument for (2). Let $[m]\in
\mcL(R,2)$, and assume $\eps_{0 *}([m])=0$. To show $Q^1_i(m) \in
\widehat H_1(\Z_2;R)$, we only have to show
$2q(v^{(1)}_i(m))=0 $ for each $i$. Set $v_i = v^{(1)}_i(m)$. This  time
set
$a_{ij}=\{1+*\}b(v_i, v_j)$, so that $2q(v_i)= a_{i
i}$. Again from \eqref{veqn1}, \eqref{veqn2}, \eqref{veqn3} we get
\begin{align*}
a_{00}&= a_{00}^2 + a_{10}^2 t\\
a_{11}&= a_{10}^2 + a_{11}^2 t .
\end{align*}
This yields again that $2q(v_1)=a_{1 1}=0, 2q(v_0)=a_{0 0}=0 $ or
$1$, and  $a_{10} = b(v_1,v_2) = 0$. But if
$2q(v_0)$ were equal to  $1$, we would recall $\eps_{0 *}([m])=0$ in
$\mcL(\Z, 2)$, and conclude
$1=
\eps_0(1)=\eps_0(Q_0^{(1)}([m])= Q_0^{(1)}(\eps_{0*}[m])=Q_0^{(1)}(0)$,
which  is impossible. This proves
$Q^{(1)}_i(\widetilde{\mcL}(R, 2))\subset \widehat H_1(\Z_2;R)$.
\end{proof}

\begin{lemma}\label{Liso}
\begin{enumerate}
\item For
$n\geq1,\;    \mcL^{\ev}(R, 2^n)\xrightarrow{\iota_n}
\mcL^{\eps}(R,2^n)$ is an isomorphism.

\item For  $n\geq
2,\quad \mcL(R, 2^{n-1})\xrightarrow{j_{n-1}}\mcL^{\ev}(R,
2^n)$ is an isomorphism.
\end{enumerate}
\end{lemma}

\begin{proof} Throughout this proof we use without comment Corollary 
\ref{colengthcor} which states that a submodule $N$ of a length one module
$M$ has colength one if and only if $N= \overline{N}$.  For example, a 
subLagrangian $S$ satisfies $S = \overline S$.

(1)\qua We construct an inverse
$G_n\co \mcL^{\eps}(R, 2^n)\to \mcL^{\ev}(R, 2^n)$  to $\iota_n$ for each
$n\geq 1$. Let $[m]=[M, b, q]\in\mcL^{\eps}(R, 2^n)$. Set $v_i
=v_i^{(n)}(m)$. By definition of $\mcL^{\eps}(R, 2^n)$ and by
Proposition \ref{onvi}(3),
$$
T(m)=\overline{\langle v_0, v_1\rangle}
$$
is a subLagrangian. (If $n=1$ and $R=\Z[t]$,  we must also use the
conclusion
$b(v_1, v_0) =0 $ from the proof of Lemma \ref{imq}). By the definition of
$v_i$, we see
$2^nq(x)=0
$ for all
$x\in T(m)^\perp$, and so $m_{T(m)}$ is of exponent $2^n$ with even type.
We define
\[
g(m)= [m_{T(m)}]\in \mcL^{\ev}(R, 2^n).
\]
Note $g(m\oplus m')=g(m)+g(m')$. Moreover if $m=(M, b, q)$
       has a Lagrangian $L$, then by Proposition \ref{onvi}(2), $T(m)\subset L$, and  so
$L/T(m)$ is a Lagrangian for $m_{T(m)}$. Therefore $g(m)=0$. Therefore $g$
induces a well defined homomorphism,
\[
G_n\co  \mcL^{\eps}(R, 2^n)\to \mcL^{\ev}(R, 2^n),\qquad  G_n([m])=
[m_{T(m)}].
\]
       Since $m$ and $ m_{T(m)}$ share the same Witt class in
$ \mcL^{\eps}(R, 2^n), \; \iota_n G_n =\id $. When $[m] \in \mcL^{\ev}(R, 2^n)$,
note $T(m)=0,\; m= m_{T(m)}$, and so
$G_n\iota_n = \id $.  Therefore $G_n$ is the inverse of $\iota_n$.

(2)\qua Suppose $ n\geq 2$. We construct an inverse $F_n\co \mcL^{\ev}(R,
2^{n})\to \mcL(R, 2^{n-1})$ to $j_{n-1}$.
Let $m=(M, b, q)\in\Lambda^{\ev}(R, 2^n)$.
Then
\[
S(m)= \overline{2^{n-1}M}
\]
        is a subLagrangian. Its subLagrangian construction
$ m_{S(m)}$ has exponent $2^{n-1}$ because $2^{n-1}(S(m)^\perp)\subset
2^{n-1}M\subset S(m)$. We  therefore define a homomorphism
\[
f\co \Lambda^{\ev}(R, 2^{n})\to\mcL (R, 2^{n-1}); \quad f(m)=
[m_{S(m)}].
\]

\begin{claima}
  Let $S$ be any  subLagrangian of $m\in \Lambda^{\ev}(R,
2^n)$ such that
\begin{equation}\label{deff}
2^{n-1}S^\perp \subset S\subset S(m).
\end{equation}
Then $f(m)=[m_S]$.  (The  first inclusion  ensures that
$m_S\in \Lambda(R,2^{n-1})$).
\end{claima}

\begin{proof}[Proof of Claim A]

A Lagrangian for $$(S^\perp/S, b_S,q_S)\oplus
(S(m)^\perp/S(m), -b_{S(m)}, -q_{S(m)})$$ is
\[
\Delta  = \{[x]\oplus [x]: x\in S(m)^\perp\}.
\]
The second inclusion of (\ref{deff}) ensures that $\Delta \subset
           S^\perp/S\oplus
S(m)^\perp/S(m).$ This
proves Claim A.
\end{proof}

\begin{claimb}
If $m\in \Lambda^{\ev}(R,
2^n)$ and $m$ has a Lagrangian, then $f(m)=0$.
\end{claimb}

\begin{proof}[Proof of Claim B]
Suppose $m=(M,b,q)$ is a quadratic linking form,  and
$2^{n }M=0,\;
          2^{n }q=0$ and $L$ is a Lagrangian  for $ m$. We must prove  $f(m)=0$.

Let $S= L\cap S(m)$, where $S(m) = \overline{2^{n-1}M}$. Obviously $S\subset S(m)$, and $L/S$ is a
Lagrangian for $m_S$. Therefore,
by Claim A, we can show $f(m)= [m_S]=0$ by showing $2^{n-1}S^\perp
\subset S$, which we now do.

          By Proposition \ref{colengtheq},     $\overline{2^{n-1}M} =
(2^{n-1}M)^{\perp\perp} = (\,_{2^{n-1}}M)^\perp$. So,
\begin{multline*}
S^\perp =(L\cap\overline{2^{n-1}M})^\perp = (L^\perp\cap
\overline{2^{n-1}M})^\perp = (L^\perp\cap ({}_{2^{n-1}}M)^\perp)^\perp \\
=
(L+{}_{2^{n-1}}M)^{\perp\perp}=\overline{L+{}_{2^{n-1}}M}.
\end{multline*}
\[\text{Therefore,\quad }
2^{n-1}S^\perp = 2^{n-1}(\overline{L+ {}_{2^{n-1}}M})\subset
\overline{2^{n-1}L}\subset L\cap \overline{2^{n-1}M}=S.
\]
This proves $2^{n-1}S^\perp \subset S$, and therefore proves Claim B.
\end{proof}

By Claim B, $f$ induces a homomorphism:
\[F_n\co \mcL^{\ev}(R, 2^{n})\to \mcL(R,
2^{n-1});\qquad F_n([m])=[m_{S(m)}].
\]
Since $m$ and $m_{S(m)}$ share   the same Witt class in $\mcL^{\ev}(R,
2^n)$, we see $j_{n-1}F_n =\id$. If $[m]\in \mcL(R, 2^{n-1})$, then
$S(m)=0$ and $m=m_{S(m)}$. So $F_n j_{n-1}=\id$.  This shows $F_n $ is
inverse  to
$j_{n-1}$.
\end{proof}

\proof[Proof of Theorem \ref{ses2q}]

By the definition and Lemmas \ref{Liso}, \ref{imq}(1),
$\mcL(R, \langle 2 \rangle)$ is the direct limit of the maps
$\iota_{n+1}\circ j_n\co \mcL(R, 2^n)\to
\mcL(R, 2^{n+1})$, and these maps are isomorphisms for $n\geq 2$. So 
$i_n\co \mcL(R, 2^n)\to \mcL(R,  \langle 2\rangle) $  is an isomorphism if
$n\geq2$. Define
\[
Q\co \mcL(R, \langle 2\rangle)\to \Z_2\;  \text{ by:}\quad Q\circ i_2 =
Q^{(2)}_0\co \mcL(R, 4)\to
\Z_2.
\]
From Lemmas \ref{Liso} and \ref{imq}, the following sequence is exact
$$
0\to \mcL(R, 2 )\xrightarrow{i_1}
\mcL(R, \langle 2 \rangle)\xrightarrow{Q}\Z_2\to 0.\eqno{\qed}
$$

\begin{proof}[Proof of Theorem \ref{devissage}]

Let  $R=\Z[t]$ or $\Z[t^-]$.  From Theorem  \ref{ses2q}
we get a commutative diagram with exact rows, whose vertical maps are
induced by the augmentation
$\eps_0$:
\[
\begin{CD}
0@>>>  \mcL(R, 2)@>i_1>>\mcL(R,  \langle 2\rangle)@>Q>>\Z_2@>>>0\\
@VVV    @VV\eps_0V
@VV\eps^{(2)}_0V@VV\id V @VVV\\
0@>>>  \mcL(\Z, 2)@>i_1>>\mcL(\Z,  \langle 2\rangle)@>Q>>\Z_2@>>>0\\
\end{CD}
\]
Since $\widetilde{\mcL}(R, 2)= \ker{\eps_0}$ and $\widetilde{\mcL}(R,  \langle
2\rangle)=
\ker{\eps^{(2)}_0}$, a simple diagram chase yields
$
i_1\co \widetilde{\mcL}(R, 2) \cong
\widetilde{\mcL}(R,  \langle 2\rangle)$.
\end{proof}

\section{Proof of Theorem \ref{computation}}\label{section-answer}

In this section we prove Theorem
\ref{computation}  by computing
$\widetilde{\mcL}(\Z[t^\pm],2)$,  the reduced Witt group of quadratic linking
forms over  $\Z[t^\pm]$ of exponent 2.

We begin with a piece of notation:
For $p(t) \in \ztt$, define $p_{\od}(t),p_{\ev}(t) \in \ztt$ by the equation
$$
p = p_{\ev}^2 + tp_{\od}^2=p_{\ev}(t^2) + tp_{\od}(t^2) \in \ztt.
$$
\begin{definition}
Let $[m]\in \ltt$. Set $v_i = v_i^{(1)}([m])$. Define
\[
B = (B_1,B_2)\co  \widetilde{\mcL}(\Z[t^\pm],2) \to 
t\Z_2[t]\times  t\Z_2[t]
\]
\begin{align}
B_1([m]) &= q(v_0)+tq(v_1);\\
B_2([m]) &=(t\,(q(v_1)_{\od}))^2 + t\,(q(v_0)_{\od})^2.
\end{align}
Note $B$ is a homomorphism. By \eqref{qnz}, $\widetilde{\mcL}^{\eps}(\zt, 2)\!=
\!\ker{(B)}.$  (Indeed, if $B_2([m])\! = 0$, then $q(v_i) = (q(v_i)_{\ev})^2$ for $i = 0,1$.  If, in 
addition, $B_1([m]) = 0$, then $q(v_0)$ and $q(v_1)$ are zero.)  By Lemma \ref{imq}, $ B$ takes values  in $\widehat
H_1(\Z_2;R)\times
\widehat H_1(\Z_2;R)t=\Z_2[t]\times t\Z_2[t]$. But since $\eps_0\circ
Q^{(1)}=Q^{(1)}\circ \eps_{0*}$, and
$\eps_{0*}([m])=0$,  $B$ takes values in $t\Z_2[t]\times t\Z_2[t]$.
\end{definition}

The following example shows  that $  \widetilde{\mcL}(\zt,
2)\xrightarrow{B}t\ztt\times t\ztt $ is an epimorphism.

\begin{example} \label{char_example}
Let $p,g\in \Z[t]$. For the quadratic linking form
$$\mcN_{p,g}  =
\left(
\ztt^2,
\begin{pmatrix}
p/2 & 1/2\\
1/2& 0
\end{pmatrix},
\begin{pmatrix}
         p/2\\
         g
\end{pmatrix}
\right)$$
we have $v_0 = (0,p_{\ev})$ and $v_1 =( 0,p_{\od})$ by equation
\eqref{veqn1}.  Hence,
\begin{align}
q(v_0)  = p_{\ev}^2g,\;\; q(v_1) &= p_{\od}^2 g.\label{qnpg}\\
B([\mcN_{p,g}]) &= (p g, (g_{\od})^2 tp). \quad \text{ In
particular},\notag
\\
       \qquad B([\mcN_{tp,1}])&= (tp,0); \\
B([\mcN_{1, tp}]-[\mcN_{t, p}]) &=(0,tp).\label{Bformula}
\end{align}
\end{example}

\proof[Proof of Theorem \ref{computation}]$\phantom{9}$

       (1)\qua The above discussion, together with Lemma \ref{Liso}(1) and Theorem \ref{even}(1) (where $\cal P$ is defined),
show that the following sequence is exact:
\[
0\to t\Z_2[t]/(V_2-1)t\Z_2[t]\xrightarrow{\iota_1\circ \mcP}
\ltt\overset{B}{\longrightarrow}t\ztt\times t\ztt\to 0.
\]

In particular, $\ltt$ is  generated by those elements  $[\mcN_{p,g}]$, for which either $p$ or
$q$ is divisible by $t$ (cf Lemma \ref{formalproperties}(3)). By Lemma
\ref{formalproperties}, parts (2), (3), and (4), 
 the map
\[
j_1\co t\Z_4[t]/(2V_2-2)t\Z_4[t] \to \ltt; \qquad j_1[tp] = [\mcN_{tp, 1}]
\]
is a well-defined homomorphism of abelian groups.  It satisfies
$
j_1[2tp] = \iota_1\mcP[tp].
$  It is clearly a $\mcV$--map.  By  
Example \ref{char_example},
 $$B_1j_1[tp] = [tp]; \quad B_2j_1[tp]=0.$$

Therefore the following sequence is exact.
\begin{equation}\label{exactck}
0\to
t\Z_4[t]/(2V_2-2)t\Z_4[t]
\overset{j_1}{\longrightarrow}\ltt\overset{B_2}{\longrightarrow}t\ztt\to 0
\end{equation}
We next claim that the map $B_2\co \ltt\to  [t\ztt]$ is a $\mcV$--map,
where  $[t\ztt]$ is the $\mcV$--module  defined
by \eqref{vstructure}.  In fact,
\begin{align*}
V_{2n+1}B_2([\mcN_{p,g}]) &= tp(t^{2n+1})\,(t^ng_{\od}(t^{2n+1}))^2
=\\B_2([\mcN_{p(t^{2n+1}), g(t^{2n+1})}]) &=
B_2V_{2n+1}([\mcN_{p, g}]);\\
V_2B_2([\mcN_{p, g}]) = 0 &= B_2([\mcN_{p(t^2),g(t^2)}])= B_2V_2([\mcN_{p,g}]).
\end{align*}
proving the claim. Therefore \eqref{exactck} becomes an exact sequence of
$\mcV$--modules.  A right splitting for it is the map
\[
j_2\co [t\ztt]\to \ltt; \qquad j_2[tp]=[\mcN_{1,tp}]-[\mcN_{t, p}].
\]

Note $j_2[tp +t2q]=[\mcN_{1,tp}]-[\mcN_{t, p}] = j_2[tp]$, so $j_2$ is well defined.  Similarly, $j_2$ is a
homomorphism because
$$
j_2[tp] + j_2[tq] -j_2[tp+tq] = [\mcP_{tp,tq}] - [\mcP_{t^2p,q}] = 0,
$$by Lemmas
\ref{formalproperties}(8) and \ref{formalproperties}(6).  It is a splitting
by
equation \eqref{Bformula}. It is a $\mcV$--map by Lemma
\ref{formalproperties}(7) and the following calculation:
\begin{align*}
V_2(j_2(tp))&=[\mcN_{1, t^2p^2}]-[\mcN_{t^2, p^2}]=0 \text{ by Lemma }
\ref{formalproperties}(5). \\
j_2 V_2 (tp) &=0 \text{, for all } tp\in [t\ztt] \text{ by }
\eqref{vstructure}.\\
V_{2n+1}(j_2(tp))&= [\mcN_{1, \, t^{2n+1}p(t^{2n+1})} ]-[\mcN_{t^{2n+1},\,
p(t^{2n+1})}]\\
&=[\mcN_{1, \, t^{2n+1}p(t^{2n+1})} ]-[\mcN_{t,\, t^{2n}
p(t^{2n+1})}] \text{ by Lemma }\ref{formalproperties}(5).\\
&=j_2V_{2n+1}(tp) \text{ by } \eqref{vstructure}.
\end{align*}
Therefore, we get  an isomorphism of
$\mcV$--modules, proving Theorem \ref{computation}(1),
\[
j_1 + j_2 \co  \frac{t\Z_4[t]}{(2V_2-2)t\Z_4[t]}
\oplus [ t\Z_2[t]] \to \widetilde \mcL (\zt,\langle 2 \rangle).
\]

(2)\qua  If $[m]\in \widetilde{\mcL}(\Z[t^-], 2),$ then $q(v^{(1)}_0(m)) \in t^2\Z_2[t^2]$,
since
$\eps_{0 }q(v^{(1)}_0(m))=0$ and  $q(v^{(1)}_0(m))\in
\widehat H_1(\Z_2; \Z[t^-])=\Z_2[t^2]$, by Lemma \ref{imq}(3). The map 
$Q^{(1)}_0\co \widetilde{\mcL}(\Z[t^-], 2)\to t^2\Z_2[t^2]$ is a homomorphism
by Lemma 
\ref{formalproperties}(2).  The map
$Q^{(1)}_0$ is injective by Theorem \ref{even}(2) and Lemma \ref{Liso}(1).
 In  fact $Q^{(1)}_0$ is an
isomorphism because $
Q^{(1)}_0([\mcN_{t^2p(t^2),1}])= t^2p(t^2),
$
by equation \eqref{qnpg}. The inverse  isomorphism to $Q^{(1)}_0$ is the
$\mcV_{\od}$--map:
$$
j\co t^2\Z_2[t^2]\to \widetilde{\mcL}(\Z[t^-], 2); \qquad j(t^2p(t^2))=
[\mcN_{t^2p(t^2), 1}].\eqno{\qed}
$$

The reader may note that $\unil{_1(\Z;\Z,\Z)}$ has not been discussed here.
In fact it is already known to be zero. See \cite{Cu}, or \cite{CK}, or
\cite{CR}. But for the
reader's convenience  we include a short proof.

\begin{prop} \label{un1}$\unil_1(\Z;\Z,\Z)=0.$
\end{prop}
\begin{proof} By Theorem \ref{unil=nl} and Proposition \ref{local},
$
\unil_1(\Z;\Z,\Z)= NL_1(\Z) = NL_1(\Z,\langle 2 \rangle),
$
which can be identified with a subgroup of 
$$
\mcL^-(\Z[t],\langle 2\rangle) = QW(\Z[t],(\Q[t]/\Z[t])^-,
\bfM\langle 2\rangle),
$$ 
the Grothendieck group of \emph{skew} quadratic linking forms over $\Z[t]$--modules $M$ with $2^nM =
0$ for some $n$.

We will show $\mcL^-(\Z[t],\langle 2\rangle)= 0$. Let $R=\Z[t]$. For any $(R,(\Q R/R)^-, \bfM\langle 2\rangle)$--form, say
$m= (M,\lambda,\mu)$, we see $\mu\co  M\to \frac{(\Q R/R)}{(1-*)(\Q
R/R)}=\frac{(\Q R/R)}{2(\Q R/R)}=0$. So  $\mu(x) =0$, and
      $\lambda(x,x)=0$ for all $x\in M$. If $2^nM=0, n\geq 2$, this implies that
$ S=\overline{2^{n-1}M}$ is a subLagrangian for $m$, and
$m_S$ has exponent $2^{n-1}$. And if $2M=0$,  it implies
that for
\emph{any} $x\neq 0$  in $M$, $S=\overline{\langle x\rangle}$ is a
subLagrangian for $m$, with $\rank_{R/2R}(m_S)\leq \rank_{R/2R}(m)-1$. Of
course this means $[m]=0$, as desired.
\end{proof}

\appendix

\section{Appendix: On length and colength and the
subLagrangian construction}\label{app}

       Here we prove several properties used in the text about the  colength  of
an $R$--submodule. Length and colength were defined in Section \ref{defs}.
We also prove that the subLagrangian construction produces a
nonsingular form in the same Witt class.

\begin{prop}\label{tripl}
Let $R$ be  a ring with involution which is
torsion free as an abelian group. Then for any $n$,  $(R, \Q R/R,
\bfM_{2^n})$ satisfies the conditions of Definition \ref{triple}, where $\bfM_{2^n}$
is the class of length one $R$--modules $M$ with $2^nM =0$.
\end{prop}

\begin{proof}
\ref{triple}(2) is obvious.

      The first part of \ref{triple}(3) is a standard fact from homological algebra.  The proof of the second part of \ref{triple}(3)  amounts to
showing that $\ext^1_R(M, \Q R/R)=0$    for any $M\in \bfM_{2^n}$. But
$\ext^i_R(M, \Q R)=0$ for all $i\geq 0$, because $2^nM=0$ and
$\times 2^n \co  \Q R\to \Q R$ is an isomorphism. Therefore $\ext^1_R(M, \Q
R/R)\cong \ext^2_R(M, R) = 0$, because $M$ has length one.

Now we prove \ref{triple}(1). Let $M\in \bfM_{2^n}$. Clearly
$2^nM^\wedge =0$. There is an exact sequence $0\to F_1\to
F_0\to M\to 0$, where
$F_0, F_1$ are finitely generated and free. Since $M^*:=\Hom_R(M,R)=0$, we
get a resulting exact sequence $0\to F_0^*\to F_1^*  \to \ext^1_R(M,
R)=M^\wedge \to 0$. This shows $M\in \bfM_{2^n}$. Let $d\co F\to F^{**}$
denote the double duality isomorphism for any finitely generated free
$R$--module. By dualizing the above short exact sequence a second time, we
get a  commutative diagram with exact rows:
\[
\begin{CD}
0@>>> F_1@>>>F_0@>>>M@>>>0\\
@VVV    @VVd_{F_0}V
@VVd_{F_0}V@VVDV @VVV\\
0@>>>  F_0^{**}@>i_1>>F_1^{**}@>>>M^{\wedge\wedge}@>>>0\\
\end{CD}
\]
where $D$ denotes the double duality map.
Since $d_{F_i}$ is an isomorphism for $i=0,1$, we conclude that $D$ is too.
\end{proof}

\begin{lemma}\label{sublag}
Let $(R, A, \bfM)$ be any  triple satisfying the conditions of
Definition \ref{triple}.  Let $S$ be a subLagrangian for
some $(R, A, \bfM)$--form
$m= (M,\lambda, \mu)$. Then $m_S$ (as defined in Section \ref{defs}) is a $(R, A, \bfM)$--form, and $[m]=[m_S]
$ in $QW(R, A, \bfM)$.
\end{lemma}

\begin{proof}  Since $M,M/S \in \bfM$ and  $0 \to S \xrightarrow{i} M
\to M/S \to
0$ is exact, we see $S \in\bf M$ and $M^\wedge \to S^\wedge$ is
surjective.  Consider the commutative diagram, where $j$ and $k$ denote inclusions:
$$
\begin{CD}
0@>>>S @>(\ad\lambda)\circ i>>M^\wedge @>k^\wedge>> (S^\perp)^\wedge  @>>>
0\\ @.   @VV jV  @VV\ad(\lambda)^{-1}V
@VV j^\wedge  V @.\\
0@>>> S^\perp @>k>> M@>i^\wedge \circ(\ad\lambda)>>S^\wedge@>>>0\\
\end{CD}
$$

The bottom row is exact, and thus, by \ref{triple}(3), so is its dual, the top row.   The exact sequence $0 \to S^\perp/S \to M/S \to
M/S^\perp
\to 0$ and the isomorphism $M/S^\perp \cong S^\wedge \in \bfM$ show
$S^\perp/S \in \bfM$ and hence
$\ker(j^\wedge ) \cong(S^\perp/S)^\wedge $.  The Snake Lemma
isomorphism $S^\perp/S = \cok(j) \cong \ker (j^\wedge) =
(S^\perp/S)^\wedge$ can be identified with $\ad\lambda_S$. So $m_S$
is nonsingular  and is therefore a $(R, A, \bfM)$--form. To see
$[m]=[m_S]$,  note $m_S\oplus (-m)$ has a Lagrangian; it is
$\{([x],x)\in (S^\perp/S)\times M: x\in S^\perp\}$.
\end{proof}

Throughout the rest of this appendix, $R$ denotes $\Z[t]$.

\begin{prop}\label{lengtheq}
Let $M$ be a finitely  generated $R$--module of exponent $2^n$ for some
$n$.  Then $M$ has length one if and only if for each $x\neq 0 $ in $M$, the
annihilator ideal satisfies $\ann_R(x)\subseteq 2R$.
\end{prop}

\begin{proof}
$\Rightarrow$: We will prove the contrapositive. Supposing that  $ x\neq 0$  in $M$ and
$\ann_R(x)\nsubseteq 2R$, we will show  that $\length_R(M)\geq2$.

Replacing $x$,
if necessary, by
$2^kx$, for some $k$,  we see that we may as well assume that $2x=0$. Since
$\ztt$ is a principal ideal domain, $\ann_R(x)=\langle 2,p\rangle$ for some $p\in
R\smallsetminus2R$.  Let
$f\co \ann_R(x)\to \Q[t]/\Z[t]$ be the $R$--map such that $f(2)=0$, and
$f(p)=\frac{1}{2}$. $f$ does not extend to an $R$--map, $R\to \Q[t]/\zt$.
Therefore the exact sequence, $0\to \ann_R(x)\to R\to Rx\to 0$ shows
that  $0\neq
\ext^{1}_R(R  x,
\Q[t]/\zt)\cong
\ext^2_R(R x, \Z[t]).$ Since $\ext^3_R(-,-)=0$, the exact sequence
$0\to Rx\to M\to M/Rx\to 0$ shows that $\ext^2_R(M,\zt)\neq 0$. This implies
that $\length_R(M)\geq 2$ as claimed.

$\Leftarrow$: Now assume for each $x\neq 0$ in $M$, that
$\ann_R(x)\subseteq 2R$. We show that $\length_R(M) =1$. The proof is by
induction on $n$. If $n=1$ the condition implies $M$ is a free $R/2R$--module
(because it is $R/2R$--torsion free). This implies $\length_R(M)=1$. Assume
the result is known for modules of exponent $2^{n-1}$. Consider the  exact
sequence
\[
0\to  \,_{2}M\to M\to M/_{2}M\to 0,
\]
where $\,_{2}M = \{x \in M : 2x = 0\}$.
For any $x\in M$ with $[x]\neq 0$ in $ M/_{2}M$ we note
$ \ann_R([x])\subseteq \ann_R(2x)\subseteq 2R$. This implies, by the
inductive hypothesis, that $\length_R (M/_{2}M)=1$, and also
$\length_R(\,_{2}M)= 1$. It follows that $\length_R(M)=1$.
\end{proof}

\begin{cor} \label{lengthcor}Let $M$ be a finitely generated $R$ module of
exponent
$2^n$ with length one. Any submodule of $M$ has length one. Also
$M^\wedge$ has length one.
\end{cor}
\begin{proof}The first claim is immediate from Proposition \ref{lengtheq}.
To see that $M^\wedge$ has length one, one applies Proposition
\ref{lengtheq},  noting that the annihilator ideal of any element of $~
\,_{2^n}(\Q[t]/\Z[t])$ (and therefore any element of $M^\wedge$) is
principal, generated by $2^k$ for some $k$.
\end{proof}

Recall the definition of the \emph{closure} of  a submodule $N$ of a
length one module $M$:
\[
\overline{N} = \{x\in M : px\in N, \text{  for some } p\in
R\smallsetminus2R\}
\]
       Clearly
$\overline{\overline{N}}=\overline{N}$. By Proposition \ref{lengtheq},
$\length_R(M/N)=1$ if and only if $N=\overline{N}$.

\begin{cor}\label{colengthcor}
Let $N$ be a submodule of a length one $R$--module $M$. Then $N$ has colength
one in $M$  if and only if $N=\overline{N}$. In particular, $\overline{N}$
has colength one.
\end{cor}

\begin{prop}\label{colengtheq}
Let $ M $ be a length one $R$--module of exponent $2^n$. Let $b\co M\times M\to
\Q[t]/\zt$ be a nonsingular linking form (relative to some involution on
$R$).  Let $N$ be a submodule of  $M$. Then $N$ has colength one
in $M$ if and only if $
       N=X^\perp$ for some $X$ in $M$. Also $\overline{N}=N^{\perp
\perp} $  for any $N$.
\end{prop}

\begin{proof}
Suppose $X$ is a submodule of $M$ for which $N=X^\perp$. Then
\[0\to N\to
M\xrightarrow{i_X^\wedge \circ (\ad b)} X^\wedge
\]
       is exact, where $i_X\co X\to M$ denotes the inclusion. But 
$\im(i_X^\wedge \circ (\ad b))$ has length one by Corollary \ref{colengthcor}.
Therefore
$N$ has colength one in $M$.

Conversely, suppose $N$ has colength one in
$M$. By Proposition \ref{tripl},  this implies
$M^\wedge\xrightarrow{i_N^\wedge} N^\wedge\to 0$ is exact.  Set
$X=\ker({i_N^\wedge \circ (\ad b)})$. This gives us an exact sequence  and its dual:
\begin{align*}
0\to
X\xrightarrow{i_X}&M
\xrightarrow{i_N^\wedge \circ (\ad b)} N^\wedge\to 0.\\
0\to
N\xrightarrow{i_N} &M
\xrightarrow{i_X^\wedge \circ (\ad b)} X^\wedge\to 0.
\end{align*}
The second sequence says that $N=X^\perp$.

Finally we show $\overline{N}=N^{\perp\perp}$.  By the definition of
$\overline{N}$, it is clear that $\overline{N}^\perp = N^\perp$. But we
have just seen that $\overline{N}=X^\perp$, for some $X\subset M$.
Moreover,
$X^{\perp\perp\perp}=X^{\perp} $ for all $X$. So
$N^{\perp\perp}=\overline{N}\,^{\perp\perp}=X^{\perp\perp\perp}=X^{\perp}=
\overline{N}$.
\end{proof}

\noindent {\bf Acknowledgement}\qua This work is partially supported by the National Science Foundation.

\end{document}

%% file: gtoutput.tex
\def\ifplaintex{\expandafter\ifx\csname documentclass\endcsname\relax}

\ifplaintex 
\else
\fi

\expandafter\ifx\csname beginpicture\endcsname\relax
\expandafter\ifx\csname documentclass\endcsname\relax
\input pictex \else
\input prepictex \input pictex \input postpictex \fi\fi

\def\gt{{\mathsurround=0pt\it $\cal G\mskip-2mu$eometry \&\ 
$\cal T\!\!$opology}}        

\def\gtp{{\mathsurround=0pt\it $\cal G\mskip-2mu$eometry \&\ 
$\cal T\!\!$opology $\cal P\!$ublications}}  

\def\lognumber#1{\def\thelognumber{#1}}
\def\volumenumber#1{\def\thevolumenumber{#1}}
\def\papernumber#1{\def\thepapernumber{#1}}
\def\volumeyear#1{\def\thevolumeyear{#1}}

\def\pagenumbers#1#2{\def\startpage{#1}\def\finishpage{#2}}
\def\published#1{\def\publishdate{#1}}
\def\proposed#1{\def\theproposer{#1}}
\def\seconded#1{\def\theseconders{#1}}
\def\received#1{\def\receiveddate{#1}}

\def\accepted#1{\def\accepteddate{#1}}

\def\asciiaddress#1{\def\theasciiaddress{#1}}
\def\asciiemail#1{\def\theasciiemail{#1}}

\long\def\asciiabstract#1{\long\def\theasciiabstract{#1}}

\let\\\par\let\thelognumber\relax
\let\thevolumenumber\relax\let\thepapernumber\relax
\let\thevolumeyear\relax\let\thesamplenumber\relax\let\startpage\relax
\let\finishpage\relax\let\publishdate\relax\let\receiveddate\relax
\let\reviseddate\relax\let\accepteddate\relax\let\theasciititle\relax
\let\theasciiauthors\relax\let\theasciiaddress\relax
\let\theasciiabstract\relax
\let\theasciiemail\relax\let\theshortauthors\relax\let\theshorttitle\relax

\long\def\maketitlep{   

\count0=\startpage

\gt\hfill      
\beginpicture
\setcoordinatesystem units <0.33truein, 0.33truein> point at 2.2 0.9
\setplotsymbol ({$\cal G$})
\plotsymbolspacing=9truept
\circulararc 315 degrees from 0 1 center at 0 0
\setplotsymbol ({$\cal T$})
\circulararc 315 degrees from 1 -1 center at 1 0
\endpicture
\break
{\small\ifx\thesamplenumber\relax 
Volume \else Sample
\fi\thevolumenumber\ (\thevolumeyear)
\startpage--\finishpage\nl
Published: \publishdate}
\vglue 0.5truein plus 0.4fil minus 0.1truein

{\parskip=0pt\leftskip 0pt plus 1fil\def\\{\par\smallskip}{\ifplaintex\large
\else\Large\fi\bf\thetitle}\par\medskip}   

\vglue 0pt plus 0.1fil 

{\parskip=0pt\leftskip 0pt plus 1fil\def\\{\par}{\sc\theauthors}
\par\medskip}

\vglue 0pt plus 0.1fil 

{\small\parskip=0pt\let\newline\\
{\leftskip 0pt plus 1fil\def\\{\par}{\sl\theaddress}\par}
\expandafter\ifx\theemail\relax    
\relax\else\vglue 5pt plus 0.02fil minus 2pt\def\\{\stdspace{\rm 
and}\stdspace} 
\cl{Email:\stdspace\tt\theemail}\fi
\ifx\theurl\relax                  
\relax\else\vglue 5pt plus 0.02fil minus 2pt\def\\{\stdspace{\rm 
and}\stdspace}
\cl{URL:\stdspace\tt\theurl}\fi\par}

\vglue 7pt plus 0.3fil minus 3pt

{\bf Abstract}
\vglue 5pt plus 0.1fil minus 2pt

\theabstract

\vglue 7pt plus 0.3fil minus 3pt

{\bf AMS Classification numbers}\quad Primary:\quad \theprimaryclass

Secondary:\quad \thesecondaryclass

\vglue 5pt plus 0.3fil minus 2pt

{\bf Keywords:}\quad \thekeywords

\vglue 10pt plus 0.5fil minus 5pt

{\small  Proposed: \theproposer\hfill Received: \receiveddate\nl
Seconded: \theseconders\hfill 
\ifx\reviseddate\relax                         
Accepted: \accepteddate                        
\else
Revised: \reviseddate                          
\fi}
\eject
}       

\let\maketitlepage\maketitlep
\let\maketitle\maketitlepage

\font\phead=cmsl9 scaled 950
\font=cmsl9 scaled 1050
\font\pnum=cmbx10 scaled 913
\font\lnum=cmbx10 
\font\pfoot=cmsl9 scaled 950
\font=cmsl9 scaled 1050
\ifplaintex
\headline{\vbox to 0pt{\vskip -4.5mm\line{\small\phead\ifnum
\count0=\startpage ISSN 1364-0380 (on line)
1465-3060 (printed) \hfill {\pnum\folio}\else\ifodd\count0\def\\{ }%
\ifx\theshorttitle\relax\thetitle\else\theshorttitle\fi\hfill{\pnum\folio}
\else\def\\{ and }{\pnum\folio}\hfill\ifx\theshortauthors\relax\theauthors
\else\theshortauthors\fi\fi\fi}\vss}}
\footline{\vbox to 0pt{\vglue 0mm\line{\small\pfoot\ifnum\count0=\startpage
\copyright\ \gtp\hfill\else
\gt, Volume \thevolumenumber\ (\thevolumeyear)\hfill\fi}\vss
}}
\else
\makeatletter
\def\@oddhead{{\small\lhead\ifnum\count0=\startpage ISSN 1364-0380 (on line)
1465-3060 (printed) \hfill {\lnum\number\count0}\else\ifodd\count0
\def\\{ }\ifx\theshorttitle\relax \thetitle \else\theshorttitle\fi\hfill
{\lnum\number\count0}\else\def\\{ and }{\lnum\number\count0}
\hfill\ifx\theshortauthors\relax 
\theauthors\else\theshortauthors\fi\fi\fi}}\def\@evenhead{\@oddhead}
\def\@oddfoot{\small\lfoot\ifnum\count0=\startpage\copyright\ \gtp\hfill\else
\gt, Volume \thevolumenumber\ (\thevolumeyear)\hfill\fi}
\def\@evenfoot{\@oddfoot}
\makeatother
\fi

\newwrite\gtoutfile
\long\gdef\makeheadfile{  
{\def\\{, }\def\s{ }
\immediate\openout\gtoutfile head.xxx
\immediate\write\gtoutfile{Proxy-for: \ifx\theasciiauthors\relax
\theauthors\else\theasciiauthors\fi\s<\ifx\theasciiemail\relax\theemail\else\theasciiemail\fi>}
\immediate\write\gtoutfile{\noexpand\\}
\immediate\write\gtoutfile{Authors: \ifx\theasciiauthors\relax
\theauthors\else\theasciiauthors\fi}
{\def\\{ }\immediate\write\gtoutfile{Title: \ifx\theasciititle\relax
\thetitle\else\theasciititle\fi}}
\immediate\write\gtoutfile{Subj-class: GT or SG or MG etc}
\immediate\write\gtoutfile{MSC-class: \theprimaryclass\ifx\thesecondaryclass\relax\else, \thesecondaryclass\fi}
\immediate\write\gtoutfile{Journal-ref: Geom. Topol. \thevolumenumber
(\thevolumeyear) \startpage-\finishpage}
\immediate\write\gtoutfile{Comments: Published by Geometry and Topology at}
\immediate\write\gtoutfile{\s\s http://www.maths.warwick.ac.uk/gt/GTVol\thevolumenumber/paper\thepapernumber.abs.html}
\immediate\write\gtoutfile{\noexpand\\}
\immediate\write\gtoutfile{}
\ifx\theasciiabstract\relax
\immediate\write\gtoutfile{\theabstract}\else
\immediate\write\gtoutfile{\theasciiabstract}\fi
\immediate\write\gtoutfile{}
\immediate\write\gtoutfile{\noexpand\\}
\immediate\write\gtoutfile{}
\immediate\closeout\gtoutfile}}  

\def\maketitlepage{\maketitlep\makeheadfile}
\let\maketitle\maketitlepage